\documentclass[11pt,twoside,a4paper]{article}
\usepackage{amsfonts,amssymb,amsmath,latexsym,epsfig,verbatim}

\newcommand{\<}{\langle}
\newcommand{\ra}{\rangle}

\newcommand{\tT}{\widetilde{T}}
\newcommand{\tp}{\widetilde{p}}
\newcommand{\tA}{\widetilde{A}}

\newcommand{\tF}{\widetilde{F}}

\def\pf{{\rm Pf}}
\def\id{\rm id}
\def\f{\overline{f}}

\newcommand{\De}{\Delta}

\newcommand{\de}{\delta}

\newcommand{\si}{\sigma}

\usepackage[latin1]{inputenc}
\usepackage{aecompl} 
\usepackage{amsmath,amsfonts,amssymb,amsthm}
\usepackage{epsfig}

\def\a{\alpha}
\def\b{\beta}
\def\d{\delta}

\def\s{\sigma}

\def\t{\tau}
\def\D{\Delta}
\newcommand{\gS}{\mathfrak{S}}

\newcommand{\cA}{\mathcal{A}}
\newcommand{\cM}{\mathcal{M}}

\newcommand{\cO}{\mathcal{O}}

\newcommand{\cC}{\mathcal{C}}

\newcommand{\cI}{\mathcal{I}}

\newcommand{\cL}{\mathcal{L}}

\newcommand{\cW}{\mathcal{W}}

\def\I{{\mathcal I}}
\def\L{{\mathcal L}}
\def\O{{\mathcal O}}

\newcommand{\R}{\mathbb{R}}
\newcommand{\N}{\mathbb{N}}
\newcommand{\C}{\mathbb{C}}

\newcommand{\Z}{\mathbb{Z}}
\newcommand{\K}{\mathbb{K}}
\newcommand{\PP}{\mathbb{P}}
\def\EE{\mathbb{E}}

\def\P{{\mathbb P}}

\newcommand{\sgn}{\operatorname{sgn}}
\newcommand{\zz}{\mathbb{Z}}
\newcommand{\rr}{\mathbb{R}}
\newcommand{\ii}{\mathcal{I}}
\newcommand{\al}{\alpha}
\newcommand{\ta}{\widetilde{\alpha}}
\newcommand{\be}{\beta}

\newcommand{\ga}{\gamma}
\newcommand{\ph}{\Phi}
\newcommand{\om}{\omega}
\newcommand{\la}{\lambda}
\newcommand{\nw}{\newtheorem}
\newcommand{\bg}{\begin}
\newcommand{\nn}{\mathbb{N}}
\newcommand{\ep}{\varepsilon}
\newcommand{\ve}{\varepsilon}

\title{Exit problems associated with affine reflection groups}
\author{Yan Doumerc and John Moriarty}
\begin{document}
\maketitle
\bibliographystyle{plain}

\begin{abstract}We obtain a formula for the distribution of the first exit time of Brownian motion from the alcove of an affine Weyl group. In most cases the formula is expressed compactly, in terms of Pfaffians. Expected exit times are derived in the type $\widetilde{A}$ case. The results extend to other Markov processes. We also give formulas for the real eigenfunctions of the Dirichlet and Neumann Laplacians on alcoves, observing that the `Hot Spots' conjecture of J. Rauch is true for alcoves.
\end{abstract}

\section{Introduction}

The distribution of the first exit time of Brownian motion from the interval (0,1) may be obtained by the reflection principle.
If $B$ is a Brownian motion, $T_{i}$ the hitting time of the level $i$,  $T_{0,1}:=T_{0}\wedge T_{1}$, and $\PP_x$ denotes the law of $B$ started at $x \in (0,1)$, then 
\begin{eqnarray}
\PP_x(T_{0,1}>t)&=&\sum_{n \in \zz}\left[ \PP_x(B_t \in 2n+(0,1))-\PP_x(
B_t \in 2n-(0,1)) \right]. \label{forz} \label{exp}
\end{eqnarray}
Using cancellation and the reflection principle, formula \eqref{exp} may be rewritten as $\PP_x(T_{0,1}>t)=\phi(x,t)$, where
\begin{eqnarray}
\phi(x,t)= \PP_{x}(T_{0}>t)+\sum_{n=1}^{\infty}(-1)^{n}[\PP_{x}(T_{-i}>t)-\PP_{x}(T_{i}>t)].\label{new}
\end{eqnarray}
More generally, suppose that $B$ is a standard Brownian motion in a real Euclidean space $V$. If $B$ is started inside the alcove $\cA$ of an affine reflection group acting on $V$, there exists an expression analogous to \eqref{exp} for the distribution of the first exit time of $B$ from $\cA$, which is given later in equation \eqref{barf}. The expression involves integration of the Gaussian kernel over the multi-dimensional alcove.
The aim of this paper is to give a formula analogous to \eqref{new}---that is, a formula for the exit probability of standard Brownian motion from the alcove, in terms of exit probabilities from simpler domains.
As an example, in the type $\tA$ case this formula involves only one-dimensional exit probabilities and can be written in terms of Pfaffians (see appendix for the definition of a Pfaffian).

To put our results in context, we state the following proposition. Let $B_1,\ldots,B_k$ be independent standard Brownian motions started at $x_1,\ldots,x_k \in \R$
and let $(\xi_n=e^{i2\pi B_n})_{n \in [k]}$ be their projections onto the circle, where $[k]:=\{1,\ldots,k\}$. Define the times of first collision 
\begin{eqnarray*}
\begin{array}{ll}T_{ij}=\inf\{t:B_i(t)=B_j(t)\} &
\widetilde{T}_{ij}=\inf\{t:\xi_i(t)=\xi_j(t)\}\\
T=\min\{T_{ij}:1 \leq i < j \leq k\} &
\widetilde{T}=\min\{\widetilde{T}_{ij}:1 \leq i < j \leq k\}. 
\end{array}
\end{eqnarray*}
Then 
$T$ is equal to the first exit time of $k$ dimensional Brownian motion started at $x=(x_1,\ldots,x_k)$ from a chamber of type 
$A_{k-1}$, and it was proved in \cite{doc} that
\begin{eqnarray}
\PP_x(T>t)=\left\{ \begin{array}{ll}\mbox{Pf}(p_{ij})_{i,j \in [k]} & (k \mbox{ even})\\
\sum_{l=1}^k(-1)^{l+1}\mbox{Pf}(p_{ij})_{i,j \in [k]\setminus\{l\}} & (k \mbox{ odd})
\end{array}\right.\label{regz}
\end{eqnarray}
where for $i < j$, $p_{ij}=\PP_x(T_{ij}>t)$
and $p_{ji}=-p_{ij}$.
As observed in \cite{MR1405497}, 
$\widetilde{T}$ is equal to the first exit time of the Brownian motion from an alcove of type 
$\widetilde{A}_{k-1}$. 
A special case of our main result gives a companion to \eqref{regz}:
\newtheorem{y1}{Proposition}
\bg{y1}\label{mano}

(i) If $k$ is even then $$\PP_x(\tT>t)=\mbox{Pf}(\tp_{ij})_{i,j \in [k]}$$ where for $i < j$,
$\widetilde{p}_{ij}=\PP_x(\widetilde{T}_{ij}>t)$ and $\widetilde{p}_{ji}=-\widetilde{p}_{ij}$. 

(ii) If $k$ is odd then 
$$\PP_x(\tT>t)=\sum_{l=1}^k(-1)^{l+1}\mbox{Pf}(q_{ij})_{i,j \in [k]\setminus\{l\}}$$
where for $i < j$,
$q_{ij}=\PP_x(\widetilde{T}_{ij}>t)+2\PP_x(\widetilde{T}_{ij}\leq t, \widetilde{T}_{ij}<T_{ij})$ and $q_{ji}=-q_{ij}$. 
\end{y1}
\noindent The relationship between $q_{ij}$ in Proposition \ref{mano} and $p_{ij}$ in \eqref{regz}
is clarified by noting that $B_{i}-B_{j}$ is also a Brownian motion and
\begin{eqnarray} \label{otter}
\psi(x,t)&:=&\PP_x(T_{0,1}>t)+2\PP_x(T_{0,1}\leq t, T_{1}<T_{0}) \nonumber\\&=&\PP_{x}(T_{0}>t)+\sum_{n=1}^{\infty}[\PP_{x}(T_{-i}>t)-\PP_{x}(T_{i}>t)],
\end{eqnarray}
which is proved in Lemma \ref{cougz} and may be compared with \eqref{new}. In the case $k=3$, $\widetilde{T}$ equals the first exit time of Brownian motion from an equilateral triangle, which is the alcove of type $\tA_2$. This relates to scaling limits occurring in, for example, a three player gambler's ruin problem and a three tower problem \cite{MR2023644,MR1975514}. As a further example, if the alcove is of type $\widetilde C_{k}$ then $\tT$ relates to first collision times for $k$ independent standard Brownian motions on the interval.

The expected exit time is obtained in the type $\tA$ case, and also a generalisation of de Bruijn's formula for multiple integrals involving determinants.
The present work extends that in \cite{doc}, where the authors consider the exit time from a chamber - that is, an unbounded domain which is the fundamental
region of a finite reflection group. The extension to $\tA_{k-1}$ with odd $k$ was prompted by Neil O'Connell, who suggested the solution for $k=3$.

The rest of the paper is organised as follows. Sections \ref{sirez} and \ref{poliz} present necessary background material, the main results with applications, a general reflection principle and an affine generalisation of De Bruijn's formula. Details of the main result in the different type cases are given in section \ref{four}. Proofs are contained in section \ref{profz}, and the real eigenfunctions of the Laplacian with Dirichlet and Neumann boundary conditions are considered in section \ref{eigz}.

\section{The geometric setting}\label{sirez}
\subsection{Finite Weyl groups and chambers}
\label{leeb}
Background on root systems and reflection groups may be found in, for example, \cite{MR1066460}.
Let $V$ be a real Euclidean space with a positive symmetric
bilinear form $\<x,y\ra$. Let $\Phi$ be an irreducible crystallographic root system in $V$
with associated reflection group $W$. Let $\D$ be a simple system in $\Phi$ 
with corresponding positive system $\Phi^+$ and fundamental chamber
$$\cC=\{ x\in V:\ \forall \; \a\in\D ,\, \<\a,x\ra>0 \}.$$ We will call $\Phi^\vee$ the set of coroots
$\a^\vee=2\a/\<\a,\a\ra$ for $\a\in \Phi$. Then, $L:=\Z-$span of $\Phi^\vee$ is a $W$-stable lattice called the coroot lattice. For $\a \in \Phi$ and $x \in V$ we make the definitions
\begin{eqnarray*}
H_{\a}&=&\{y \in V: \< \a,y \ra =0\} \\s_{\a}(x)&=&x-\<\al,x \ra \a^{\vee}.
\end{eqnarray*}
Thus $s_{\a}$, $\a\in\Phi^+$ are the reflections in
$W$. 
\subsection{Affine Weyl groups and alcoves} \label{affw}
The affine Weyl group $W_a$ asociated with $\Phi$ is the group generated by
all affine reflections with respect to the hyperplanes 
$H_{(\a,n)}=\{x\in V :\
\<x,\a\ra=n\}$, $\a\in \Phi^{+}$, $n\in\Z$. It has a semi-direct
product decomposition in terms of the Weyl group $W$ and the coroot lattice $L$: each element of $W_a$ may be written uniquely as
$\t(l) w$, where $w\in W$ and $\t(l)$ is
the translation by $l\in L$. We may therefore attribute a sign to each $w_a=\t(l) w\in W_a$ by $\varepsilon(w_a)=\varepsilon(w):=\det(w)$. The fundamental alcove is the bounded domain defined by 
\begin{eqnarray*}
\cA&=&\{x\in V :\ \forall \;
\a\in \Phi^+,\, 0<\<x,\a\ra<1\}\\ &=&\{x\in V : \<x,\widetilde{\alpha}\ra<1\text{ and } \forall \;
\a\in \D,\, \<x,\a\ra >0\}
\end{eqnarray*}
where $\widetilde{\alpha}$ is the highest positive
root. 

\subsection{Affine root systems}
We refer to \cite{MR1838580} for this formalism although we use slightly modified
notations for the sake of consistency.
\nw{y2}[y1]{Definition}
\begin{y2}
If $\Phi$ is an irreducible crystallographic root
system as previously introduced, the corresponding affine root system is
$\Phi_a:=\Phi\times \Z$. For $\la=(\a,n)\in\Phi_a$ and $x \in V$ we
define 
\begin{eqnarray*}
\la(x)=\la.x&=&\<\a,x\ra-n \\
H_{\la}&=&\{y \in V: \la.y=0\} \\
s_\la(x)&=&x-(\la.x)\a^\vee
\end{eqnarray*}
\end{y2}
Thus $s_{\la}$ is the reflection with respect to the hyperplane $H_\la$, and we may write $s_{\la}=\tau(n \a^\vee)s_{\al}$.
Writing $w_a=\tau(l)w \in W_a$, we have that $W_a$ acts on $V$ by $w_a(x)=w(x)+l$; we define further the action of $W_a$ on 
$\Phi_a$ by 
\nw{vegc}[y1]{Definition}
\bg{vegc} For $w_a=\tau(l)w \in W_a$ and $\la=(\a,n)\in\Phi_a$,
$$w_a(\la)=(w\a,n+\<w\a,l\ra)\in\Phi_a.$$
\end{vegc}
We then have $w_a(\la) . x=\la . w_a^{-1}(x)$ for $w_a\in W_a,\,
\la\in\Phi_a,\, x\in V$, which is analogous to the fact that $W$ is a group of isometries; we also have $w_a H_\la=H_{w_a(\la)}$. 
If $\la=(\a,m),\, \mu=(\beta,n)\in \Phi_a$ then we will refer to the angle between $\la$ and $\mu$, meaning the angle between $\a$ and $\b$; by $\la\perp\mu$ we mean $\<\alpha,\beta\ra=0$. The usual properties of a reflection are then preserved:
$s_\la(\la)=(-\a, -n)=:-\la$ and $s_\la(\mu)=\mu$ if $\la\perp\mu$.

\nw{yan2}[y1]{Definition}
\begin{yan2}
The affine simple system is $\D_a:=\{(\a,0),\, \a\in\D;\ (-\widetilde{\al},-1) \}$ and the
corresponding positive system is $\Phi_a^{+} :=\{(\a,n):\ (n=0
\text{ and } \a\in \Phi^+) \text{ or } n\le-1\}$. 
\end{yan2}

This definition is tailor-made so that $$\cA=\{x\in V :\ \forall \
\la\in \Phi_a^{+},\, \la(x)>0\}=\{x\in V :\ \forall \ \la\in \D_a,\, \la(x)>0\}.$$

\label{trag}

\section{Background and main results}\label{poliz}
We present here our main results, which extend the main result in \cite{doc} to the affine cases.
In section \ref{corz} we give some applications
in the type $\tA$ case.\subsection{Consistency}
Let $(\mathcal{W}, \phi, \phi^+,\delta,F) \in \{(W,\Phi,\Phi^+,\Delta,\cC), (W_a,\Phi_a,\Phi^+_a,\Delta_a,\cA)\}$ and 
for $I \subset \phi^+$
define $\mathcal{W}^I=\{w \in \mathcal{W}:wI \subset \phi^+ \}$ and $\mathcal{I}=\{wI:w \in \mathcal{W}^I\}$.
For $S \subset \phi$, we define 
the set of orthogonal subsets of $S$:
$$\O(S):=\{Y\subset S:\ \forall \; \la \not=\mu \in Y,\, \la\perp \mu\}.$$
\nw{var}[y1]{Definition}
\bg{var}[\textbf{Consistency}]
\label{varz}
\begin{itemize}
\item We will say that $I$ satisfies
hypothesis (C1) if there exists $J \in \mathcal{O}(\delta \cap I)$ such that
if $J \subset A \in \ii$ then $A=I$.
\item We will say that $I $ satisfies
hypothesis (C2) if the restriction of the determinant to the subgroup $U=\{w \in \mathcal{W}: wI=I\}$ is trivial, i.e. $\forall \; w \in U, \ep(w)= \det w =1$.
\item We will say that $I$ satisfies hypothesis (C3) if $\ii$ is finite.
\item $I$ will be called \textbf{consistent} if it satisfies (C1), (C2) and
(C3).
\end{itemize}
\end{var}
Condition (C2) makes it
possible to attribute a sign to every element of $\I$ by $\varepsilon_A:=\varepsilon(w)$ 
for $A \in \I$, where $w$ is any element of $W^I$ with $w I=A$.

\subsection{Reflectability}
Let $X=(X_t,t \geq 0)$ be a $V$-valued process and let $\PP_x$ denote the law of $X$ started at $x \in F$. 
We will call $X$ \textbf{reflectable} if it satisfies the conditions of the following:
\nw{wvar}[y1]{Definition}
\bg{wvar}[\textbf{Reflectable process}]
\label{wvarz}
\begin{itemize}
\item $X$ has the strong Markov property.
\item The sample paths of $X$ are almost surely continuous.
\item The law of $X$ is $\mathcal{W}$-invariant - that is,
$\PP_x \circ (w X)^{-1}=\PP_{w x} \circ X^{-1}$ for all $w \in \mathcal{W}, x \in V$.
\end{itemize}
\end{wvar}
\subsection{Exit times} \label{notoo}
We now introduce some notation for exit times.
Let $X$ be a reflectable process in $V$. For convenience we may 
write each $\lambda \in \phi^+$ in the form $(\al,n)$ by 
identifying $\al \in \Phi$ with $(\al,0) \in \Phi_a$. Then for $\lambda=(\al,n) \in \phi^+$ 
define $T_{\la}=\inf\{t \geq 0: \la.X_t=0\}$ and for $A=\{\la_{1},\ldots,\la_{k}\} \subset \phi^+$ write 
$T_{A}:=T_{\la_{1},\ldots,\la_{k}}:=\min_{\la \in A}T_{\la}$.
Finally, let $T$ denote the first exit time of $X$ from the fundamental chamber $\cC$---that is, $T=T_{\d}$ in the finite case $(\mathcal{W}, \phi^+,\delta,F)=(W,\Phi^+,\Delta,\cC)$; and let $\tT$ denote the first exit time of $X$ from the fundamental alcove $\cA$---that is, $\widetilde{T}=T_{\delta}$ in the affine case $(\mathcal{W}, \phi^+,\delta,F)= (W_a,\Phi^+_a,\Delta_a,\cA)$.

\subsection{Main results} \label{dropz}
The following Theorem extends the main result of \cite{doc} to include those affine Weyl groups which have a consistent subset; the details of its application to particular affine Weyl groups are given in section \ref{four}. Theorem \ref{var3z} deals with an important case where a consistent subset is not available.
\nw{var2}[y1]{Theorem}
\bg{var2}
\label{var2z} 
Suppose $I$ is consistent, $X$ is reflectable and $x \in F$. Then :
\begin{eqnarray}\PP_x(T_{\delta}>t)=
\sum_{A \in \ii}\ep_A \PP_x(T_A>t). \label{cannnon}
\end{eqnarray}
\end{var2}
Note that the sum is finite even for
affine Weyl groups. In the $\widetilde{A}_{k-1}$ case with odd $k$, no consistent subset is available and we require a different formalism: for $A \in \mathcal{O}(\ph)$, define
\begin{eqnarray*}
E_A&=&\{v\in \text{Span}(A) : \forall \; \beta \in A,\, (v,\beta)\in\Z\}
\\
\ep^A_v&=& (-1)^{\#\{\be \in A \; : \; \<v,\be\ra>0\}}\\
|v|_A&=&\max\{|\<v,\beta\ra|:\beta \in A\}
\end{eqnarray*}
where $\#$ is the cardinality function. For $v,\be \in V$ define 
\begin{eqnarray*} \label{previ}
T_{\be,v}=\inf 
\{t \geq 0:\<X_t,\be\ra=\<v,\be\ra\}, \qquad
T_{A,v}=\min_{\be \in A}T_{\be,v}.
\end{eqnarray*}
To clarify, $E_A$ is a lattice (equal to the $\Z$-span of $A/2$) and $\ep^A_{\cdot}$, $|\cdot|_A$ give a sign and norm respectively on this lattice; and $T_{A,v}$ is the first time that the projections of $X_{t}$ and $v$ coincide along  some $\be \in A$.
\nw{var3}[y1]{Theorem}
\bg{var3}
\label{var3z} 
In the case $\mathcal{W}=\widetilde{A}_{k-1}$ with $k$ odd,
if $X$ is reflectable and $x \in \cA$ then
\begin{eqnarray}
\PP_x(\widetilde{T}>t)=\sum_{A \in \ii}\sum_{k \in \nn}\sum_{\substack{v \in E_A \\ |v|_A=k}}\ep_A \ep^A_v \PP_x(T_{A,v}>t)\label{hjp}
\end{eqnarray}
if this sum converges, where $I$ and $\ii$ are taken from the case $\cW=A_{k-1}$.
\end{var3}

\subsubsection{The `orthogonal' case} \label{orthogz}
We begin this section by recording some definitions.
\nw{semio}[y1]{Definition} 
\bg{semio}
\begin{itemize}
\item We say $A \subset \phi^{+}$ is \textbf{block-orthogonal} if it can be partitioned into 
blocks $(\rho_{i})$ such that $\rho_{i} \perp \rho_{j}$ for $i\neq j$ and each $\rho_{i}$ is either a singlet or a pair of 
roots whose mutual angle is $\pi$.
\item We say $A \subset \phi^{+}$ is \textbf{semi-orthogonal} if it can be partitioned into 
blocks $(\rho_{i})$ such that $\rho_{i} \perp \rho_{j}$ for $i\neq j$ and each $\rho_{i}$ is either a singlet or a set of 
vectors whose mutual angles are integer multiples of $\pi/4$.
\end{itemize}
\end{semio}
If $I$ is block-orthogonal and $X$ has independent components in orthogonal directions, \eqref{cannnon} factorises to give
\begin{eqnarray}\PP_x(T_{\delta}>t)=
\sum_{A \in \ii}\ep_A \prod_{i}\PP_x(T_{\rho_{i}}>t). \label{facz}
\end{eqnarray}
In many cases it is convenient to write \eqref{facz} in terms of Pfaffians, and the details are given in section \ref{four}. Under slightly stronger conditions on $X$, \eqref{hjp} factorises analogously:
\nw{damo}[y1]{Proposition}
\bg{damo}
\label{ninenine}
In the case $\cW=\tA_{k-1}$ with $k$ odd, under conditions on $X$ which hold for Brownian motion we have
\begin{equation}\label{leen}\PP_x(\tT>t) =\sum_{A \in \ii}\ep_A \prod_{\be \in A}
\left(\PP_x[T_{\be}\wedge T_{(\be,1)}>t]+2\PP_x[T_{\be}>T_{(\be,1)}\leq t]
\right).
\end{equation}
\end{damo}
This expression may also be written in terms of Pfaffians, as noted in Proposition \ref{mano}(ii).

\subsection{Applications}
\label{corz}

\subsubsection{Expected exit time in the type $\tA$ case}
The fundamental chamber for $A_{k-1}$ is $\cC=\{x \in V:x_1>x_2>\ldots>x_k \}$ where $V=\rr^k$ or $V=\{x \in \rr^k: x_1+\ldots+x_k=0\}.$ As noted in the introduction, $T$ is the first `collision time' between any two coordinates of $X$.
The fundamental alcove for the corresponding affine Weyl group
$\widetilde{A}_{k-1}$ is $\cA=\{x \in V:1+x_k>x_1>x_2>\ldots>x_k \}$.

In the $A_{k-1}$ case, an explicit formula for the expected exit time of Brownian motion from the fundamental chamber has 
been obtained in \cite{doc}:
\begin{eqnarray} \label{vanm}
\EE_x(T)=\sum_{\pi \in P_2(k)}(-1)^{c(\pi)}F_{p}(x_{\pi})
\end{eqnarray}
where $p=\lfloor k/2 \rfloor$ and $x_{\pi}=(x_i-x_j)_{\{i<j\}\in \pi} \in \rr^p_+$. Here $P_2(k)$ is the set of partitions of $[k]=\{1,\ldots,k\}$ into $k/2$ pairs if $k$ is even and into $(k-1)/2$ pairs and a singlet if $k$ is odd. The quantity $c(\pi)$ is the number of crossings in the partition $\pi$ (if k is odd, we consider an extra pair made of the singlet and another singlet labelled 0, and use this pair to compute the number of crossings); for an illustration see section \ref{drop}.
The notation $\{i<j\}\in \pi$ means that $\{i,j\}\in \pi$ and $i<j$,
 and the function $F_{p}$ is given by
$$F_p(y_1,\ldots,y_p)= 
\frac{2^{p+1}\Gamma(p/2)}{\pi^{p/2}(p-2)}
\int_0^{y_1}\ldots \int_0^{y_p}\frac{dz_1\ldots dz_p}
{(z_1^2+\ldots+z_p^2)^{p/2-1}}.
$$
We prove an analogous formula:

\nw{libr}[y1]{Proposition} 
\bg{libr}In the $\widetilde{A}_{k-1}$ case, if $X$ is Brownian motion then $$\EE_x(\tT)=\sum_{\pi \in P_2(k)}(-1)^{c(\pi)}\tF_{p}(x_{\pi})$$ where
\label{librz}
$$\tF_p(y_1,\ldots,y_p)=
\frac{2^{2p}}{\pi^{p+2}}
\sum_{l \in \mathbb{O}^p}\frac{1}{(l_1^2+\ldots+l_p^2)}\prod_{s=1}^p \frac{1}{l_s}\sin (\pi l_s y_s)
$$
where $\mathbb{O}=2\nn+1$ if $k$ is even and 
$\mathbb{O}=2\nn$ if $k$ is odd, and 
with the definition 
$\frac{1}{l_s} \sin (\pi l_s y_s) = \frac{1}{2} \pi y_s$ when $l_s=0.$

\end{libr}
In the case $k=3$ we will recover the known formula \begin{eqnarray}
\mathbb{E}_x(\tT)=x_{12}x_{23}(1-x_{13}),\label{alwayz}
\end{eqnarray}
where $0<x_{ij}=x_i-x_j<1$,
for the expected exit time of Brownian motion from an equilateral triangle.

\subsubsection{Dual formulae and small time behaviour}

Dual to \eqref{new} and \eqref{otter} are the formulae
\begin{eqnarray}
1-\phi(x,t)&=& \PP_{x}(T_{0}\leq t)+\sum_{n=1}^{\infty}(-1)^{n}[\PP_{x}(T_{-i} \leq t)-\PP_{x}(T_{i}\leq t)] \\
1-\psi(x,t)&=&\PP_{x}(T_{0}\leq t)+\sum_{n=1}^{\infty}[\PP_{x}(T_{-i}\leq t)-\PP_{x}(T_{i}\leq t)].
\end{eqnarray}
In the block-orthogonal case of section \ref{orthogz}, these dual formulae may be used to obtain asymptotics for the small time behaviour of the exit probability. For example, exact asymptotics can be obtained in the Brownian case, as in section 4.6.2 of \cite{doc} (we omit the details).

\subsubsection{Eigenfunctions for alcoves}

In section \ref{eigz}, using results from \cite{MR570879}, we obtain formulae for the real eigenfunctions of the Laplacian on alcoves with Dirichlet or Neumann boundary conditions. This confirms a version of the `Hot Spots' conjecture of J. Rauch for alcoves. We also prove the following

\nw{yan6}[y1]{Proposition}
\begin{yan6} \label{trigoharmoyou}
Let $\cA$ be the fundamental alcove of an affine Weyl group, and let the corresponding Weyl group have positive system $\Phi^+$.
The function 
\begin{equation}\label{definitrigo}
H(x):=\prod_{\a \in \Phi^+} \sin\left(\pi \langle x,\a\rangle \right)
\end{equation} is an eigenfunction for the Laplacian with Dirichlet boundary conditions on $\cA$.
Since $H$ is positive on $\cA$, it is the principal eigenfunction. Further, each eigenfunction is divisible by $H$ in the ring of trigonometric polynomials.
\end{yan6}

\subsection{The reflection principle and De Bruijn Formula}

In this section we recall a reflection principle in the context of finite or affine reflection groups, and use it to deduce a generalisation of 
a formula of De Bruijn for evaluating multiple integrals involving determinants. For the proof of Theorem \ref{loiz} we refer to \cite{MR1678525} and references therein.

\nw{loi}[y1]{Theorem}
\bg{loi}
\label{loiz}
Let $\PP_x$ denote the law of a reflectable process $X$ started
from $x \in F$. Then for all measurable sets $B \subset F,$
\begin{eqnarray}
\PP_x[X_t \in B, T_{\de}>t]= \sum_{\om \in \mathcal{W}} \ep(\om)
\PP_x[X_t \in \om B]. \label{barf}
\end{eqnarray}

\end{loi}

We apply this result in the following propositions, whose applications include the evaluation of Selberg type integrals of eigenfunctions of the Dirichlet Laplacian on an alcove (see section \ref{eigz}). 

Suppose $I$ is consistent. For $A \in \mathcal{I},$ denote by $W_A$ the group generated by the reflections $s_{\la}, \;\la \in A.$ Denote by $F_A$ the fundamental region associated
to $A$, $F_A=\{x \in V: \forall \; \la \in A, \la(x)>0\}$. Also, since $\Phi=\Phi^+ \cup (-\Phi^+)$, for $\be \in \Phi$ and $B \subset \Phi$ we may define the absolute values
\begin{eqnarray}
|\be|= \left\{ \begin{array}{ll} 
\be & : \be \in \Phi^+ \\ -\be &: -\be \in \Phi^+ \end{array}
\right., \qquad |B|=\{|\be|: \be \in B\}. \label{sweeney}
\end{eqnarray}

Assume that $F_A$ is the fundamental region for the reflection group $W_A$, which is certainly the case if $I$ is block-orthogonal or semi-orthogonal.
Theorems \ref{var2z} and \ref{loiz} in the Brownian case 
give

\nw{cert}[y1]{Proposition}

\bg{cert}
\label{certz}
If $I$ is consistent and $f:V \to \rr$ is integrable, then 
\begin{eqnarray}
\label{funn} \int_F \sum_{w \in \cW} \ep(w)f(wy)dy=\sum_{A \in \ii}\ep_A \sum_{w \in W_A}\ep(w)
\int_{F_A}f(wy)dy.
\end{eqnarray}
\end{cert}

In many cases, if $f$ factorises 
this formula may be expressed in terms of Pfaffians (see \cite{doc}); the type $A$ case
was first obtained by de Bruijn~\cite{MR0079647} using different methods. The next two results work out the corresponding results in the type $\widetilde{A}$ case.

\nw{green}[y1]{Proposition}
\bg{green}
\label{greenz} 
Let $\mathcal{W}=\widetilde{A}_{k-1}$ and let $f(y_1,\ldots,y_k)
=f_1(y_1)\ldots f_k(y_k)$ for integrable functions 
$f_i:\rr \to \rr$. If $k$ is even then
\begin{eqnarray*}
\int_{\cA} \sum_{\om \in W_a}\ep(\om)f(\om y)dy=
\pf(J_{ij})_{i,j \in [k]}
\end{eqnarray*}
where 
$J_{ij}= \int (-1)^{\lfloor y-z \rfloor}f_i(y) f_j(z) dy dz.$
\end{green}

\nw{green2}[y1]{Proposition}
\bg{green2}
\label{green2z} 
Under the conditions of Proposition \ref{greenz}, if $k$ is odd then
\begin{eqnarray*}
\int_{\cA} \sum_{\om \in W_a}\ep(\om)f(\om y)dy=
\sum_{l=1}^k(-1)^{l+1}\int_{\rr}f_l \pf(H_{ij})_{i,j \in [k]\setminus \{l\}}
\end{eqnarray*}
if
$$\sum_{m=1}^{\infty}\int_{y-z \in (-\infty,-m)\cup(m,\infty)}|f_i(y)f_j(z)|dy dz < \infty,$$
where
$H_{ij}= \int \sgn (y-z) f_i(y) f_j(z) dy dz + 2 \sum_{m=1}^{\infty}
\int_{y-z \in (-\infty,-m)\cup(m,\infty)} \sgn (y-z) f_i(y) f_j(z) dy dz.$
\end{green2}

\section{Application to the different type cases} \label{four}
Throughout this section we will assume that $X$ has independent components in orthogonal directions, to enable the writing of formula \eqref{cannnon} in terms of Pfaffians.
\subsection{The $\widetilde{A}_{k-1}$ case, $k$ even} \label{drop}
In this case, $W$ is $\mathfrak{S}_k$ acting on $\R^k$ by permutation
of the canonical basis vectors, $V=\R^k$ or $V=\{x\in\R^k :\ \sum_i
x_i=0\}$, $\Phi^+=\{e_i-e_j, 1 \leq i < j \leq k\}$, $\D=\{e_i-e_{i+1},\,1\le i\le k-1\}$,
$\widetilde{\alpha}=e_1-e_k$, $\cA=\{x\in V :\ 1+x_k>x_1>\cdots>x_k\}$, 
$\a^\vee=\a$ for $\a\in \Phi$ and $L=\{d\in\Z^k :\
\sum_{i=1}^k d_i=0\}$.\\
For even $k=2p$, we take $I=\{ (e_{2i-1}-e_{2i},0),\,
(-e_{2i-1}+e_{2i},-1)\, ;\, 1\le i\le p\}$. Then $I$ is consistent and block-orthogonal, and
$\I$ can be identified with the set $P_2(k)$ of
partitions of $[k]$ as shown in the following example for
$k=4$. Under this identification, the sign $\varepsilon_A$ is just the parity of
the number $c(\pi)$ of crossings.

\unitlength=.8cm
\begin{figure}[h]
\begin{picture}(11,14.5)
\put(4,13){\circle*{.15}}
\put(5,13){\circle*{.15}}
\put(6,13){\circle*{.15}}
\put(7,13){\circle*{.15}}
\put(3.85,12.5){1}
\put(4.85,12.5){2}
\put(5.85,12.5){3}
\put(6.85,12.5){4}
\put(3.5,11.8){$\pi=\{\{1,4\},\{2,3\}\}$}
\put(3.5,11.1){$A=\{(e_1-e_4,0),(e_2-e_3,0),$}
\put(3.5,10.4){$(-e_1+e_4,-1),(-e_2+e_3,-1)\}$}
\put(3.5,9.7){$c(\pi)=0$}
\qbezier(4,13)(5.5,14)(7,13)
\qbezier(5,13)(5.5,13.5)(6,13)

\put(4,8){\circle*{.15}}
\put(5,8){\circle*{.15}}
\put(6,8){\circle*{.15}}
\put(7,8){\circle*{.15}}
\put(3.85,7.5){1}
\put(4.85,7.5){2}
\put(5.85,7.5){3}
\put(6.85,7.5){4}
\put(3.5,6.8){$\pi=\{\{1,3\},\{2,4\}\}$}
\put(3.5,6.1){$A=\{(e_1-e_3,0),(e_2-e_4,0),$}
\put(3.5,5.4){$(-e_1+e_3,-1),(-e_2+e_4,-1)\}$}
\put(3.5,4.7){$c(\pi)=1$}
\qbezier(4,8)(5,9)(6,8)
\qbezier(5,8)(6,9)(7,8)

\put(4,3){\circle*{.15}}
\put(5,3){\circle*{.15}}
\put(6,3){\circle*{.15}}
\put(7,3){\circle*{.15}}
\put(3.85,2.5){1}
\put(4.85,2.5){2}
\put(5.85,2.5){3}
\put(6.85,2.5){4}
\put(3.5,1.8){$\pi=\{\{1,2\},\{3,4\}\}$}
\put(3.5,1.1){$A=\{(e_1-e_2,0),(e_3-e_4,0),$}
\put(3.5,.4){$(-e_1+e_2,-1),(-e_3+e_4,-1)\}$}
\put(3.5,-.3){$c(\pi)=0$}
\qbezier(4,3)(4.5,4)(5,3)
\qbezier(6,3)(6.5,4)(7,3)
\end{picture}
\caption{Pair partitions and their signs for $\widetilde{A}_3$.}
\label{figA}
\end{figure}

Hence, the formula can be written as
\begin{equation}\label{forma}
\P_x(T>t)=\sum_{\pi\in P_2(k)} (-1)^{c(\pi)} \prod_{\{i<j\}\in\pi} \widetilde{p}_{ij}={\rm Pf}\, \left(\widetilde{p}_{ij}\right)_{i,j\in [k]}
\end{equation}
where
$\widetilde{p}_{ij}=\P_x(T_{(e_i-e_j,0),(-e_i+e_j,-1)}>t)=\P_x(\forall
s\le t,\, 0<X^i_s-X^j_s<1)=\phi(x_i-x_j,2t)$ where $\phi(x,t)$ is defined in \eqref{fzt}.\\
\\
For odd $k$, we do not have a consistent subset as
the sign $\varepsilon_A$ is
not well-defined. The difference between even and odd $k$ can
be seen directly at the level of pair partitions: interchanging $1$ and $k$ in the
blocks of $\pi\in P_2(k)$ (which corresponds to the reflection with respect to $\{x_1-x_k=1\}$, which is the affine hyperplane of the alcove) changes the sign of $\pi$ if $k$ is even
while the sign is unaffected if $k$ is odd. In this case
(which includes, for example, the equilateral triangle in the case $\widetilde{A}_2$), we instead use Theorem  \ref{var3z}.

\subsection{The $\widetilde{C}_{k}$ case}
In this case, $W$ is the group of signed permutations acting on $V=\R^k$, $\D=\{2e_k,e_i-e_{i+1},\,1\le i\le k-1\}$,
$\widetilde{\alpha}=2e_1$, $\cA=\{x\in \R^k :\ 1/2>x_1>\cdots>x_k>0\}$ and $L=\Z^k$.\\
For even $k=2p$, we take $$I=\{ (e_{2i-1}-e_{2i},0),\,
(2e_{2i},0),\,(-2e_{2i-1},-1) ;\, 1\le i\le p\}.$$ 
For odd $k=2p+1$, $$I=\{ (e_{2i-1}-e_{2i},0),\,
(2e_{2i},0),\, (-2e_{2i-1},-1),\,(2e_k,0),\, (-2e_k,-1) ;\, 1\le i\le p\}.$$
$I$ is semi-orthogonal and again, $\I$ can be identified with $P_2(k)$;  the formula is
\begin{equation}\label{formc}
\P_x(T>t)=\sum_{\pi\in P_2(k)} (-1)^{c(\pi)}\, \check{p}_{s(\pi)} \prod_{\{i<j\}\in\pi} \check{p}_{ij}
\end{equation}
where
$$\check{p}_{ij}=\P_x(T_{(e_i-e_j,0),(2e_j,0),(-2e_i,-1)}>t)=\P_x(\forall
s\le t,\, 1/2>X^i_s>X^j_s>0),$$ $$\check{p}_{i}=\P_x(T_{(2e_i,0),(-2e_i,-1)}>t)=\P_x(\forall
s\le t,\, 1/2>X^i_s>0),$$ 
and $s(\pi)$ is the singlet of $\pi$, the term $\check{p}_{s(\pi)}$ being absent for even $k$.\\
\\
Everything can be rewritten in terms of Pfaffians:
\begin{equation}\label{pfac}
\P_x(T>t)=\left\lbrace
\begin{array}{cc}
{\rm Pf}\, \left(\check{p}_{ij}\right)_{i,j\in [k]} & \text{ if $k$ is even},\\
\sum_{l=1}^k (-1)^{l-1}\,\check{p}_{l}\,{\rm Pf}\,\left(\check{p}_{ij}\right)_{i,j\in
[k]\backslash \{l\}}& \text{ if $k$ is odd}.
\end{array} \right.
\end{equation}

\nw{yan5}{Remark}
\begin{yan5}
This formula can be obtained directly by applying the exit probability formula for
the chamber of type $C_k$ (which is the same as $B_k$) to the Brownian motion killed when reaching $1/2$. But it was natural to include it in our framework.
\end{yan5}

\subsection{The $\widetilde{B}_{k}$ case}
$W$ is the group of signed permutations acting on $V=\R^k$, $\D=\{e_k,e_i-e_{i+1},\,1\le i\le k-1\}$,
$\widetilde{\alpha}=e_1+e_2$, $\cA=\{x\in \R^k :\ x_1>\cdots>x_k>0,\ x_1+x_2<1\}$ and
$L=\{d\in\Z^k :\ \sum_i d_i \text{ is even}\}$.\\
For even $k=2p$, we take $$I=\{ (e_{2i-1}-e_{2i},0),\,
(e_{2i},0),\,(-e_{2i-1}-e_{2i},-1) ;\, 1\le i\le p\}.$$ 
For odd $k=2p+1$, $$I=\{ (e_{2i-1}-e_{2i},0),\,
(e_{2i},0),\,(-e_{2i-1}-e_{2i},-1),\,(e_k,0),\,(-e_k,-1) ;\, 1\le i\le p\}.$$
In this case, $I$ is semi-orthogonal and the formula is:
\begin{equation}\label{formb}
\P_x(T>t)=\sum_{\pi\in P_2(k)} (-1)^{c(\pi)} \bar{p}_{s(\pi)}\, \prod_{\{i<j\}\in\pi} \bar{p}_{ij}
\end{equation}
where
\begin{eqnarray*}
\bar{p}_{ij}&=&\P_x(T_{(e_i-e_j,0),(-e_i-e_j,-1),(e_j,0)}>t)=\P_x(\forall
s\le t,\, 1-X_s^j>X^i_s>X^j_s>0),\\
\bar{p}_{i}&=&\P_x(T_{(e_i,0),(-e_i,-1)}>t)=\P_x(\forall
s\le t,\, 1>X^i_s>0)
\end{eqnarray*}
and $s(\pi)$ denotes the singlet of $\pi$, the
term $\bar{p}_{s(\pi)}$ being absent for even $k$.\\
\\
Everything can be rewritten in terms of Pfaffians:
\begin{equation}\label{pfab}
\P_x(T>t)=\left\lbrace
\begin{array}{cc}
{\rm Pf}\, \left(\bar{p}_{ij}\right)_{i,j\in [k]} & \text{ if $k$ is even},\\
\sum_{l=1}^k (-1)^{l-1}\,\bar{p}_l\,{\rm Pf}\,\left(\bar{p}_{ij}\right)_{i,j\in
[k]\backslash \{l\}}& \text{ if $k$ is odd}.
\end{array} \right.
\end{equation}

\subsection{The $\widetilde{D}_{k}$ case}
$W$ is the group of evenly signed permutations acting on $V=\R^k$,
$\D=\{e_i-e_{i+1},e_{k-1}+e_k,\ 1\le i\le k-1\}$,
$\widetilde{\alpha}=e_1+e_2$, $\cA=\{x\in \R^k :\ x_1>\cdots>x_{k-1}>|x_k|,\ x_1+x_2<1\}$ and
$L=\{d\in\Z^k :\ \sum_i d_i \text{ is even }\}$.\\
For even $k=2p$, we take $$I=\{ (e_{2i-1}-e_{2i},0),\,(-e_{2i-1}+e_{2i},-1),\,(e_{2i-1}+e_{2i},0),\,(-e_{2i-1}-e_{2i},-1) ;\, 1\le i\le p\}.$$ 
For odd $k=2p+1$, $$I=\{ (e_{2i}-e_{2i+1},0),\,(-e_{2i}+e_{2i+1},-1),\,(e_{2i}+e_{2i+1},0),\,(-e_{2i}-e_{2i+1},-1) ;\, 1\le i\le p\}.$$ 
$I$ is block-orthogonal and the formula then becomes:
\begin{equation}\label{formd}
\P_x(T>t)=\sum_{\pi\in P_2(k)} (-1)^{c(\pi)} \prod_{\{i<j\}\in\pi} \breve{p}_{ij}
\end{equation}
where
\begin{eqnarray*}
\breve{p}_{ij}&=&\P_x(T_{(e_i-e_j,0),(-e_i+e_j,-1),(e_i+e_j,0),(-e_i-e_j,-1)}>t)=\acute{p}_{ij}\,\grave{p}_{ij},\\
\acute{p}_{ij}&=&\P_x(\forall
s\le t,\, 1>X^i_s-X^j_s>0)=\phi(x_i-x_j,2t),\\
\grave{p}_{ij}&=&\P_x(\forall
s\le t,\, 1>X^i_s+X^j_s>0)=\phi(x_i+x_j,2t)
\end{eqnarray*}
and $\phi(x,t)$ is defined in \eqref{fzt}. Everything can be rewritten in terms of Pfaffians:
\begin{equation}\label{pfad}
\P_x(T>t)=\left\lbrace
\begin{array}{cc}
{\rm Pf}\, \left(\breve{p}_{ij}\right)_{i,j\in [k]} & \text{ if $k$ is even},\\
\sum_{l=1}^k (-1)^{l-1}\,{\rm Pf}\,\left(\breve{p}_{ij}\right)_{i,j\in
[k]\backslash \{l\}}& \text{ if $k$ is odd}.
\end{array} \right.
\end{equation}

\subsection{The $\widetilde{G}_2$ case}
Here, $V=\{x\in \R^3, \sum_i x_i=0\}$,
$\Phi^+=\{e_3-e_1,e_3-e_2,e_1-e_2,-2e_1+e_2+e_3,-2e_2+e_1+e_3,2e_3-e_1-e_2\}$,
$\widetilde{\alpha}=2e_3-e_1-e_2$, $\D=\{e_1-e_2,-2e_1+e_2+e_3\}$ and
$L=\{d\in V: \forall i,\,3d_i\in \Z\}$.\\
We take
$I=\{(e_1-e_2,0),(-e_1+e_2,-1),(2e_3-e_1-e_2,0),(-2e_3+e_1+e_2,-1)\}$,
which is consistent and we can describe $\I$ as $\{I,A_1,A_2\}$ with
$A_1=\{(e_3-e_1,0),(-e_3+e_1,-1),(-2e_2+e_1+e_3,0),(2e_2-e_1-e_3,-1)\}$,
$\varepsilon_{A_1}=-1$, $A_2=\{(e_3-e_2,0),(-e_3+e_2,-1),(-2e_1+e_2+e_3,0),(2e_1-e_2-e_3,-1)\}$,
$\varepsilon_{A_2}=1$. 

In this case, the chamber $\cA$ is a triangle $ABC$ with angles $(\pi/2,\pi/3,\pi/6)$ as represented in Figure \ref{figGz}. If $T_R$ denotes the exit time from the region $R$ of the plane and $\P(R)=\P_x(T_{R}>t)$, then Theorem \ref{var2z} in this case gives 
\begin{equation}\label{trirec}
\P(ABC)=\P(ADEC)-\P(FJCG)+\P(FHCI),
\end{equation}
where $ADEC$, $FJCG$, $FHCI$ are rectangles, as shown in Figure 11.2.
\unitlength=1cm
\begin{figure}[h]
\begin{picture}(10,5)
\put(7,0){\line(1,0){3.46}}
\put(7,0){\line(-1,0){3.3}}
\put(7,2){\line(1,0){3.46}}
\put(7,2){\line(-1,0){3.3}}
\put(7,4){\line(1,0){3.46}}
\put(7,4){\line(-1,0){3.3}}
\put(3.7,0){\line(0,1){4}}
\put(7,0){\line(0,1){4}}
\put(10.46,0){\line(0,1){4}}
\qbezier(7,4)(7,4)(9.309,0)
\qbezier(7,0)(7,0)(9.309,4)
\qbezier(9.309,0)(9.309,0)(10.46,2)
\qbezier(9.309,4)(9.309,4)(10.46,2)
\qbezier(7,4)(7,4)(10.46,2)
\qbezier(7,0)(7,0)(10.46,2)
\qbezier(7,4)(7,4)(4.961,0)
\qbezier(7,0)(7,0)(4.961,4)
\qbezier(4.961,0)(4.961,0)(3.7,2)
\qbezier(4.961,4)(4.961,4)(3.7,2)
\qbezier(7,4)(7,4)(3.7,2)
\qbezier(7,0)(7,0)(3.7,2)
\put(6.65,1.65){$A$}\put(8.35,1.65){$B$}\put(6.85,4.1){$F$}
\put(7,2){\circle*{.1}}\put(8.15,2){\circle*{.1}}\put(7,4){\circle*{.1}}
\put(10.52,2){$D$}\put(10.52,0){$E$}\put(6.85,-.4){$C$}
\put(10.46,2){\circle*{.1}}\put(10.46,0){\circle*{.1}}\put(7,0){\circle*{.1}}
\put(5.45,3.08){\circle*{.1}}\put(5,3){$G$}
\put(5.45,.92){\circle*{.1}}\put(5,.7){$I$}
\put(8.75,3.03){\circle*{.1}}\put(8.85,3){$H$}
\put(8.75,1){\circle*{.1}}\put(8.85,.8){$J$}
\end{picture}
\caption{Tiling associated with $\widetilde{G}_2$}
\label{figGz}
\end{figure}
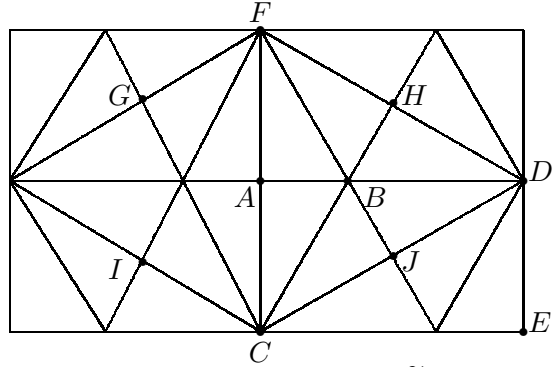

\subsection{The $\widetilde{F}_4$ case}
Recall that $V=\R^4$, $\Phi^+=\{e_i\pm e_j,\, 1\le i<j\le 4;\,e_i
,\, 1\le i\le
4;\,(e_1\pm e_2\pm e_3\pm e_4)/2\}$,
$\D=\{e_2-e_3,e_3-e_4,e_4,(e_1-e_2-e_3-e_4)/2\}$, $\widetilde{\alpha}=e_1+e_2$ and
$L=\{d\in\Z^4 :\ \sum_i d_i\text{ is even}\}$.\\
$I:=\{(e_2-e_3,0),(-e_2+e_3,-1),(e_1-e_4,0),(-e_1+e_4,-1),(e_3,0),(e_4,0)\}$
turns out to be consistent and so Theorem \ref{var2z} applies, although in this case it does not seem easy to give the formula in a compact way.

\section{Proofs} \label{profz}
\subsection{Theorem \ref{var2z}}
All the formalism of affine root systems has been set for the proofs in this section to be the same as those in \cite{doc}. Therefore, we only state the lemmas (without proofs) to show how they have to be modified in this context. 
\nw{yan8}[y1]{Lemma}
\bg{yan8}
\label{bijectional}
If $I$ is consistent then for $K\subset I$ and $\la\in\delta\cap K^\perp$ we have $s_{\la} \L=\L$, where
$$\L=\{ A\in \I:\ K\subset A,\ \la\notin A\} .$$
\end{yan8}
\nw{yan9}[y1]{Lemma}
\begin{yan9}\label{ggal}
Suppose condition (C3) is satisfied and that the function $f:\I\to\R$ and the root $\la\in\delta$ are such that
$f(A)=0$ whenever $\la\in A$, and
$f(A)=f(s_{\la} A)$ whenever $\la\notin A$. Then
$\sum_{A\in \I} \varepsilon_A f(A) =0.$
\end{yan9}
\nw{yan10}[y1]{Lemma}
\begin{yan10}\label{gggal}
If $I$ is consistent then we have:
$\sum_{A\in \I} \varepsilon_A =1.$
\end{yan10}

\textbf{{Proof of Theorem \ref{var2z}}.} 
Appealing to Lemma \ref{gggal} and the fact that  
$T_{\delta} \leq T_A$ for all $A \in \ii$, it is equivalent to prove $\sum_{A \in \ii}\ep_A \PP_x(T_A>t,T_{\delta} \leq t)=0$ and therefore sufficient to prove $\sum_{A \in \ii}\ep_A \PP_x(T_A>t,T_{\delta}=T_{\la} \leq t)=0$
for each $\la \in \delta$. Since $X$ is reflectable, $f(A)=\PP_x(
T_A>t,T_{\delta}=T_{\la} \leq t)$ satisfies the conditions 
of Lemma \ref{ggal}.
\hfill $\Box$\\

\subsection{Theorem \ref{var3z}}
\label{delayz}
Before proving Theorem \ref{var3z} we record some preliminary results.
Since a consistent subset is available in the setting of the finite reflection group $A_{k-1}$,
we work in this context. The definitions of $V,\D,\ta$ and $\Phi^+$ when $W=A_{k-1}$ have been given in section \ref{drop}.
It is proved in \cite{doc} that $I=\{e_1-e_2,e_3-e_4,\ldots,e_{k-2}-e_{k-1}\}$ is consistent and orthogonal. In the following we will make use of the notation introduced in 
sections \ref{notoo}-\ref{dropz} and \eqref{sweeney}. Also, for $\be \in \Phi^+$ define \[\cL_{\be}=\{A \in \mathcal{I}: \be \notin A\}.\]

\nw{1tild}[y1]{Lemma}
\bg{1tild}
\label{1tildz}
$A\mapsto |s_{\ta}A |$ is a permutation of $\cL_{\ta}$ and $\ep_{|s_{\ta} A|}=(-1)^{|A \setminus \ta^\perp|+1}\ep_A$ for all $A \in \cL_{\ta}$.
\end{1tild}
\textbf{Proof}
Take $A=\om I \in \cL_{\ta}$. Since the elements of $A\setminus \ta^\perp$ are orthogonal to each other so are those of $s_{\ta}(A\setminus \ta^\perp)$ thus the product $p:=\prod_{\be \in A\setminus \ta^\perp} s_{s_{\ta} (\be)}$ is well-defined (and commutative). First, take $\al\in A\cap \ta^\perp$. Then $|s_{\ta} \al|=|\al|=\al$. If $\be \in A\setminus \ta^\perp$ then $\be\not= \al$ hence $\be\perp \al$. Together with $\ta\perp \al$, we get $s_{\ta}(\be)\perp \al$ and $s_{s_{\ta} (\be)}s_{\ta} \al=s_{s_{\ta} (\be)} \al=\al=|s_{\ta}\al|$. Thus, $p\, s_{\ta} \al=|s_{\ta}\al|$. Second, take $\al\in A\setminus \ta^\perp$. Then $s_{\ta} \al\in -\Phi^+$ and $|s_{\ta}\al|=s_{s_{\ta}(\al)}s_{\ta}\al$. For $\be\in A\setminus \ta^\perp$ and $\be\not= \al$ we have $\be\perp -\al$ so $s_{\ta}\be\perp -s_{\ta} \al=s_{s_{\ta}\al}s_{\ta}\al$. Therefore $p\, s_{\ta} \al=s_{s_{\ta}\al}s_{\ta}\al=|s_{\ta}\al|$. We have proved that $|s_{\ta}A |=p\, s_{\ta} A=p\, s_{\ta} \om I$. Together with $|s_{\ta}A| \subset \Phi^+$, this yields $|s_{\ta}A| \in \cI$ and $\ep_{|s_{\ta} A|}=(-1)^{|A \setminus \ta^\perp|+1}\ep_A$. Moreover $\ta\notin |s_{\ta}A|$ since $\ta\notin A$. Consequently  $|s_{\ta}A|\in \cL_{\ta}$. It remains to observe that $A\mapsto |s_{\ta}A |$ is an involution hence a bijection.
\hfill $\Box$\\

Observing that $\<e_1-e_k,e_1-e_j\ra=\<e_1-e_k,e_i-e_k\ra=1$ for $1<i,j<k$ gives 
\nw{reff}[y1]{Lemma}
\bg{reff}
\label{reffz}For all $\be \in \Phi^+\setminus (\ta \cup \ta^\perp)$ we have $\<\ta,\be\ra=1$.
\end{reff}

Also, calculations such as 
\begin{eqnarray*}
\ep^A_{s_{\al}v}&=&(-1)^{ \#\{\be \in A \; : \; \<s_{\al}v,\be\ra > 0\}}=
(-1)^{ \#\{\be \in s_{\al}A \; : \; \<s_{\al}v,s_{\al}\be\ra > 0\}}\\
&=&(-1)^{ \#\{\be \in s_{\al}A \; : \; \<v,\be\ra > 0\}}=\ep^{s_{\al}A}_{v} 
\end{eqnarray*}
establish
\nw{coro}[y1]{Lemma}
\bg{coro}\label{coroz}For all $\a \in \Phi,$ $A \in \cO(\Phi)$ and $v \in E_{A}$ we have $$s_{\a}E_A=E_{s_{\a}A},\  \ep^A_{s_{\a}v}=\ep^{s_{\a}A}_v\ \text{and}\ |s_{\a}v|_A=|v|_{s_{\a}A}.$$
\end{coro}

\nw{fur}[y1]{Proposition}
\bg{fur} \label{furz} 

(i) $\sum_{A \in \mathcal{I}}\sum_{k \in \nn} \sum_{\substack{v \in E_A \\ |v|_A=k}}\ep_A\ep_v^A=1$. 

(ii) Suppose $f:\ii \times V \to \rr$ is such that 
$f(A,v)=f\left(|s_{\ta} A|,p_{s_{\ta} A}(s_{\ta,1}v)\right)$ 
whenever $\ta \notin A$ ($p_B$ is the orthogonal projection on $\text{Span}(B)$) and $f$ is sufficiently decreasing in the second variable (see the precise condition (\ref{rapdec}) in the proof). Then $\sum_{A \in \cL_{\ta}}\sum_{k \in \nn} \sum_{\substack{v \in E_A \\ |v|_A=k}}\ep_A\ep_v^Af(A,v)$ converges and its sum is zero.

(iii) If $f:\ii \times V \to \rr$ and $\al \in \Delta$ are such that
$f(A,v)=f(s_{\al}A,s_{\al}v)$ whenever $\al \notin A$, then $\sum_{A \in \cL_{\al}}\sum_{k \in \nn} \sum_{\substack{v \in E_A \\ |v|_A=k}}\ep_A\ep_v^Af(A,v)$ converges and its sum is zero.
\end{fur}

\textbf{Proof} (i) 
For $A\in \cO(\Phi)$ and $\al\in A$ define $$S(A,k)=\sum_{\substack{v \in E_A\\|v|_A=k}}\ep_v^A \quad {\rm and } \quad S'(A,\a,k)=\sum_{\substack{v \in E_A\\|v|_A=k}}\mathbf{1}_{v\notin \al^{\perp}}\,\ep_v^A,$$ where $\mathbf{1}$ is the indicator function. 
Since $\ep_0^A=1$ for all $A\in\cI$ and $\sum_{A \in \mathcal{I}}\ep_A=1$ by Lemma \ref{gggal}, we have
\begin{equation}\label{yesnu} \sum_{A \in \mathcal{I}}\sum_{k \in \nn}\sum_{\substack{v \in E_A\\|v|_A=k}}\ep_A\ep_v^A =1+
\sum_{A \in \mathcal{I}}\ep_A \sum_{k \geq 1} S(A,k). \end{equation}
If $u\notin \al^{\perp}$, $\ep_{s_{\al}u}^A=\ep_u^{s_{\al}A}=\ep_u^{A\setminus \{\al\} \cup \{-\al\}}=-\ep_u^{A}$. Thus, setting $v=s_{\al}u$ in $S'(A,\a,k)$ and using $s_{\al}E_A=E_{s_{\al}A}=E_A$, $|s_{\al}u|_A=|u|_{s_{\al}A}=|u|_A$, $\mathbf{1}_{s_{\al}u\notin \al^{\perp}}=\mathbf{1}_{u\notin \al^{\perp}}$, we get
$S'(A,\a,k)=\sum_{\substack{u \in E_A\\|u|_A=k}} \mathbf{1}_{u\notin \al^{\perp}}\,\ep_{s_{\al}u}^A= -S'(A,\a,k)=0.$
Therefore $$S(A,k)=\sum_{\substack{v \in E_A \cap \al^{\perp}\\|v|_A=k}}\,\ep_v^A=\sum_{\substack{v \in E_{A\setminus\{\al\}}\\|v|_{A\setminus\{\al\}}=k}}\,\ep_v^A= S(A\setminus\{\al\},k).$$ By iteration $S(A,k)=S(\emptyset,k)$, which is an empty sum (since $E_{\emptyset}=\{0\}$ and $k\ge 1$) hence null.  

(ii) Take $A \in \cL_{\ta}$ and $u\in E_A$. Now
$\ep^{|s_{\ta}A|}_{s_{\ta,1}u}=(-1)^{ \#\{\be \in |s_{\ta}A| \; : \;
\<s_{\ta,1}u,\be\ra > 0\}}$ and $s_{\ta,1}u=s_{\ta}u+\ta$;
therefore if $\be \in |s_{\ta}A| \setminus \ta^{\perp}$ then 
 writing $\ga=-s_{\ta}\be \in A \setminus \ta^{\perp}$ and applying Lemma \ref{reffz} we have $\left(\<s_{\ta,1}u,\be\ra > 0 \iff \<u,\ga\ra <1\right)$.
Also, if $\b \in |s_{\ta}A| \cap \ta^{\perp} = A \cap \ta^{\perp}$ then
$\left(\<s_{\ta,1}u,\be\ra>0 \iff  \<u,\be\ra > 0\right)$. We conclude that
\begin{eqnarray*}
\ep^{|s_{\ta}A|}_{s_{\ta,1}u}&=& (-1)^{ \#\{\ga \in A \setminus \ta^{\perp} \; : \; \<u,\ga\ra < 1\}+\#\{\be \in A \cap \ta^{\perp} \; : \;
\<u,\be\ra > 0\}} \\ 
&=& (-1)^{ \#\{\ga \in A \setminus \ta^{\perp} \; : \; \<u,\ga\ra < 1\}+\#\{\be \in A \setminus \ta^{\perp} \; : \;
\<u,\be\ra > 0\}}\, \ep_u^A = (-1)^{|A \setminus \ta^\perp|}\,\ep_{u}^A.
\end{eqnarray*}
Since $\ep_{|s_{\ta} A|}=(-1)^{|A \setminus \ta^\perp|+1}
\ep_A$ by Lemma \ref{1tildz}, we have 
\begin{equation}\label{propintyan}
	\ep_{|s_{\ta} A|}\ep_{s_{\ta,1}u}^{|s_{\ta} A|}f(|s_{\ta} A|,p_{s_{\ta}A}(s_{\ta,1}u))=-\ep_A\ep_{u}^A f(A,u).
\end{equation}
For $K\in\N = \{0,1,2,\ldots\}$, set $$S_K=\sum_{A \in \cL_{\ta}}\sum_{k=0}^K \sum_{\substack{v \in E_A \\ |v|_A=k}}\ep_A \ep_v^A f(A,v).$$ Using the permutation $A\mapsto |s_{\ta}A|$ of $\cL_{\ta}$ from Lemma \ref{1tildz} and since both $E_{|B|}=E_B$ and $|v|_{|B|}=|v|_B$ for $B \subset \cO(\Phi)$, we get
$$S_K=\sum_{A \in \cL_{\ta}}\sum_{k=0}^K \sum_{\substack{v \in E_{s_{\ta}A} \\ |v|_{s_{\ta}A}=k }}\ep_{|s_{\ta}A|} \ep_v^{|s_{\ta}A|} f(|s_{\ta}A|,v).$$
For $A \in \cL_{\ta}$ and $u\in E_A$, define $g_A(u)=p_{s_{\ta}A}(s_{\ta,1}u)=s_{\ta}u+p_{s_{\ta}A}(\ta)$. Then $g_A(u)\in \text{Span}( s_{\ta}A)$ and for all $\be\in A$, 
$$\left\langle g_A(u),s_{\ta}\be \right\rangle=\left\langle s_{\ta,1}u,s_{\ta}\be \right\rangle=\left\langle u,\be \right\rangle-\left\langle \ta,\be \right\rangle \in\Z$$
since $u\in E_A$ and $\< \ta,\be\ra \in\{0,1\}$ (Lemma \ref{reffz}). This proves that $g_A(u)\in E_{s_{\ta}A}$ and $|g_A(u)|_{s_{\ta}A}=|u|_A+\eta_A(u)$ where $\eta_A(u)\in\{-1,0,1\}$. Then, $g_A :  E_A\to E_{s_{\ta}A}$ is easily seen to be a bijection (check that $g_A^{-1}(v)=p_A(s_{\ta,1}v)$). Using this bijection as well as (\ref{propintyan}), we obtain
\begin{equation}\label{somyan}S_K=-\sum_{A \in \cL_{\ta}}\sum_{k=0}^K \sum_{\substack{u \in E_A \\ |u|_A+\eta_A(u)=k}}\ep_A \ep_u^A f(A,u).\end{equation}
Now, for $i\in\{-1,0,1\}$, let $S_i(k)=\sum_{A \in \cL_{\ta}} \sum_{\substack{u \in E_A \\ |u|_A=k,\; \eta_A(u)=i}}\ep_A \ep_u^A f(A,u)$.  Then (\ref{somyan}) reads
$$
S_K=-\sum_{k=0}^K \big(S_0(k)+S_1(k-1)+S_{-1}(k+1)\big).$$
Since $S_K=\sum_{k=0}^K \big(S_0(k)+S_1(k)+S_{-1}(k)\big)$ by definition,
we get
$$2S_K=-S_1(-1)+S_1(K)+S_{-1}(0)-S_{-1}(K+1).$$
Now, $S_1(-1)$ and $S_{-1}(0)$ are empty sums hence null. The requirement on $f$ is 
\begin{equation}\label{rapdec}
\lim_{k\to\infty}\sum_{\substack{A\in\cI,u\in E_A\\|u|_A=k}} |f(A,u)|=0,	
\end{equation}
which clearly implies $\lim_{K\to \infty} S_i(K)=0$ and consequently $\lim_{K\to \infty} S_K=0$.\\

(iii) Since $s_{\alpha}$ is a permutation of $\cL_{\al}$ (Lemma \ref{bijectional}), 
\begin{eqnarray*}
U_K&:=&\sum_{A \in \cL_{\al}}\sum_{k=0}^K \sum_{\substack{v \in E_A \\ |v|_A=k}}\ep_A\ep_v^A f(A,v)
=\sum_{A \in \cL_{\al}}\sum_{k=0}^K \sum_{\substack{v \in E_{s_{\al}A} \\ |v|_{s_{\al}A}=k}}\ep_{s_{\al}A}\ep_v^{s_{\al}A} f(s_{\al}A,v)\\
&=&\sum_{A \in \cL_{\al}}\sum_{k=0}^K \sum_{\substack{u \in s_{\al}E_{s_{\al}A} \\ |s_{\al}u|_{s_{\al}A}=k}}\ep_{s_{\al}A}\ep_{s_{\al}u}^{s_{\al}A} f(s_{\al}A,s_{\al}u)=-U_K,
\end{eqnarray*}
where the third equality follows from setting $u=s_{\al}v$ and the fourth follows from Lemma \ref{coroz}, the given property of $f$ and $\ep_{s_{\al}A}=-\ep_A$. Thus, all partial sums $U_K$ are zero. 
\hfill $\Box$\\


\textbf{Proof of Theorem \ref{var3z}}

From (i) of Proposition \ref{furz}, the theorem is equivalent to 
$$\sum_{A \in \mathcal{I}}\sum_{k \in \nn} \sum_{\substack{v \in E_A \\ |v|_A=k}}\,\ep_A\,\ep_v^A\, \left(\PP_x[T_{A,v}>t]-\PP_x[\widetilde{T}>t]\right)=0.$$
For $A\in\cI$, $v\in E_A$ and $\be \in A$ we have $\<v,\be\ra\in\Z$ hence $\<v,\be\ra\notin (0,1)$. Thus, $\widetilde{T}\le T_{\be,v}$ and so $\widetilde{T}\le T_{A,v}$. This implies \begin{eqnarray}
\PP_x[T_{A,v}>t]-\PP_x[\widetilde{T}>t]&=&\PP_x[T_{A,v}>t,\, \widetilde{T}\le t]\nonumber\\&=&\sum_{\la \in \Delta_a} \PP_x[T_{A,v}>t,\, \widetilde{T}=T_\la\le t].\label{dag}
\end{eqnarray}
(If the events in \eqref{dag} are not disjoint (up to a set of probability 
zero) we may easily redefine the $T_{\la}$ to make them disjoint, without affecting the following reflection argument.)
Now fix $\la=(\al,n) \in \left\{\Delta\times\{0\}\right\} \cup \{(\ta,1)\}$ (this set is more convenient than  $\Delta_a$ since we have $(\ta,1)$ instead of $(-\ta,-1)$). We will prove that
$$\sum_{A \in \mathcal{I}}\sum_{k \in \nn} \sum_{\substack{v \in E_A \\ |v|_A=k}}\,\ep_A\,\ep_v^A\, \PP_x[T_{A,v}>t,\, \widetilde{T}=T_\la\le t]=0.$$ Since $\PP_x[T_{A,v}>t,\, \widetilde{T}=T_\la\le t]=\PP_x[\widetilde{T}=T_\la\le t]-\PP_x[T_{A,v}\le t,\, \widetilde{T}=T_\la\le t]$ and using (i) of Proposition \ref{furz} again, this is equivalent to $$S:=\sum_{A \in \mathcal{I}}\sum_{k \in \nn} \sum_{\substack{v \in E_A \\ |v|_A=k}}\,\ep_A\,\ep_v^A\, f(A,v)=\PP_x[\widetilde{T}=T_\la\le t],$$
where $f(A,v)=\PP_x[T_{A,v}\le t,\, \widetilde{T}=T_{\la}\le t]$. 
We first prove that 
\begin{equation}\label{invar_f}
	f(A,v)=f(s_{\al}A,s_{\la}v).
\end{equation} Since $f(A,v)=\PP_x[\widetilde{T}=T_\la\le t]-g(A,v)$ where $g(A,v)=\PP_x[T_{A,v}> t,\, \widetilde{T}=T_\la\le t]$, it is enough to prove $g(A,v)=g(s_{\al}A,s_{\la}v)$. We define $\widehat{X}_u=X_u1_{u \leq T_{\la}}+s_{\la}
X_u1_{u > T_{\la}}$ and use obvious `hat notations' for stopping times associated with $\widehat{X}$. The reflectable process $X$ has the same law as $\widehat{X}$ so that
$g(A,v)=\PP_x[\widehat{T}_{A,v}> t,\, \widehat{\widetilde{T}}=\widehat{T}_{\la}\le t]$. Since $X$ and $\widehat{X}$ coincide before $T_{\la}=\widehat{T}_{\la}$, we have $\widehat{\widetilde{T}}=\widetilde{T}$. Together with $\widehat{T}_{A,v}\,\mathbf{1}_{\widehat{T}_{A,v}>T_{\la}}=T_{s_{\al}A,s_{\la}v}\,\mathbf{1}_{T_{s_{\al}A,s_{\la}v}>T_{\la}}$, this yields
	$$g(A,v)=\PP_x[T_{s_{\al}A,s_{\la}v}>t,\, \widetilde{T}=T_{\la}\le t]=g(s_{\al}A,s_{\la}v),$$
which proves the claim.

In addition to the equality $f(A,v)=f(|A|,p_A v)$, equation (\ref{invar_f}) ensures that $f$ has the relevant property for Proposition \ref{furz} to yield
$$\sum_{A \in \cL_{\al}}\sum_{k \in \nn} \sum_{\substack{v \in E_A \\ |v|_A=k}}\,\ep_A\,\ep_v^A\,f(A,v)=0$$
so that $S=\sum_{\substack{A \in \mathcal{I}\\ \al \in A}}\ \sum_{k \in \nn} \sum_{\substack{v \in E_A \\ |v|_A=k}}\,\ep_A\,\ep_v^A\,f(A,v).$
If $\al \in A$ then  $f(A,v)=f(A,s_{\la}v)$ (thanks to (\ref{invar_f})) and if $\la(v)\neq 0$ then $\ep^A_{v}=-\ep^A_{s_{\la}v}$. Then as in the proof of Proposition \ref{furz}(ii) we can use the bijection $v\mapsto s_\la v$ to remove cancelling pairs and appeal to 
property (\ref{rapdec}) to conclude that $$\sum_{k \in \nn} \sum_{\substack{v \in E_A \\ |v|_A=k}}\,\mathbf{1}_{\la(v)\neq 0}
\ep_v^A f(A,v)=0$$ so that $S=\sum_{\substack{A \in \mathcal{I}\\ \al \in A}}\ep_A S(A)$, where $S(A):=\sum_{k \in \nn} \sum_{\substack{v \in E_A \\ |v|_A=k}}\,\ep_v^A\,\mathbf{1}_{\la(v)=0} \,f(A,v).$
If $\al \in A$ and $\la(v)=0$ then $f(A,v)=\PP_x[\widetilde{T}=T_{\la} \leq t]$ and $$S(A)=\PP_x[\widetilde{T}=T_{\la} \leq t] \sum_{k \in \nn} \sum_{\substack{v \in E_A \\ |v|_A=k}}\,\ep_v^A\,\mathbf{1}_{\la(v)=0}.$$ For $\be\in A\backslash\{\al\}$, the bijection $v\mapsto s_\be v$ flips the sign $\ep_v^A$ creating pair cancellations for the terms with $v$ not orthogonal to $\be$. Repeating this for all $\be\not=\al$ as in the proof of Proposition \ref{furz}(i), we are left only with that $v$ which is a multiple of $\al$ such that $\la(v)=0$, i.e. $v=n \al/2$ : we have $$S(A)=\ep^A_{n \al / 2}\PP_x[\widetilde{T}=T_{\la}\leq t].$$ 
It remains only to show that $\sum_{\substack{A \in \mathcal{I}\\ \al \in A}}\ep_A\ep^A_{n \al / 2}=1$. When $\al \in \Delta$ this follows from the proof of Lemma \ref{gggal}, which can be found in \cite{doc}; for $\al=\ta$, observe that
$\ep^A_{\ta / 2}=-1$ if $\ta \in A$. Identifying $A \in \ii$ with $\pi \in P_2(k)$ as in section \ref{drop}, we have $\left( \ta \in A \iff \{1,k\} \in \pi \right)$. 
Now $\{1,k\}$ crosses the pair containing 0, and no other pair. 
It follows that $c(\pi)=1+c(\pi \setminus \{1,k\})$, so 
$$\sum_{\substack{A \in \ii\\ \ta \in A}}\ep_A=
\sum_{\substack{\pi \in P_2(k)\\ \{1,k\} \in \pi}}(-1)^{c(\pi)}=
-\sum_{\pi \in P_2(k-2)}(-1)^{c(\pi)}=-1$$
by Lemma \ref{gggal}. \hfill $\Box$\\
\label{helz}

\subsection{Proposition \ref{ninenine}} 

\nw{snee}[y1]{Lemma}
\bg{snee}
\label{sneez}
For $A \in \mathcal{O}(\Phi^+)$, 
if projections of $X$ in orthogonal directions are independent then for  $x \in \cA$, 
\begin{equation}
	\sum_{k \in \nn} \sum_{\substack{v \in E_A \\ |v|_A=k}} \ep^A_v \PP_x[T_{A,v}>t]=\prod_{\be \in A}
\sum_{n \in \nn}\sum_{\substack{k\in\Z \\ |k|=n}} \si(k) \PP_x[T_{(\be,k)}>t],
\label{chin}
\end{equation}
where $\si(k)=-1$ if $k>0$ and $\si(k)=1$ otherwise,
if these sums converge.
\end{snee}
\textbf{Proof}
Set $A=\{\be_1,\ldots,\be_p\}$. Rewriting and expanding the respective partial sums gives, for $N \in \nn,$
\begin{eqnarray*}
	\prod_{i=1}^p \sum_{n=0}^N\sum_{\substack{k\in\Z \\ |k|=n}} \si(k) \PP_x[T_{(\be_i,k)}>t]
=\sum_{n=0}^N\ \sum_{\substack{\vec{k}=(k_1,\ldots,k_p)\in\Z^p\\|\vec{k}|_{\infty}=n}}\ \prod_{i=1}^p \si(k_i) \PP_x[T_{(\be_i,k_i)}>t].	
\end{eqnarray*}
Now, $\vec{k}=(k_1,\ldots,k_p)\mapsto v=\frac{1}{2} \sum_{i=1}^p k_i\be_i$ is a bijection from $\Z^p$ to $E_A$ satisfying $\<v,\be_i\ra=k_i$ so that $T_{(\be_i,k_i)}=T_{\be_i,v}$, $|v|_A=|\vec{k}|_{\infty}$ and $\ep_v^A=\prod_{i=1}^p \si(k_i)$. By independence $\prod_{i=1}^p \PP_x[T_{(\be_i,k_i)}>t]=\PP_x[\min_{i} T_{(\be_i,k_i)}>t]=\PP_x[T_{A,v}>t]$, and letting $N \to \infty$ concludes the proof. \hfill $\Box$\\

\nw{coug}[y1]{Lemma}
\bg{coug}
If $X$ is reflectable 
then for $x \in \cA$, 
$$\PP_x[T_{\be}\wedge T_{(\be,1)}>t]+2\PP_x[T_{\be}>T_{(\be,1)} \leq t] = \sum_{n \in \nn}\sum_{\substack{k\in\Z \\ |k|=n}} \si(k) \PP_x[T_{(\be,k)}>t].$$
\label{cougz}
\end{coug}
\textbf{Proof} 
Let
\begin{eqnarray*}
S1&=& \sum_{k=1}^{\infty}\left(\PP_x[T_{(\be,-k)}>t,T_{\be}\wedge T_{(\be,1)}>t]-\PP_x[T_{(\be,k)}>t,T_{\be}\wedge T_{(\be,1)}>t]\right)
\\
S2&=& \sum_{k=1}^{\infty}\left(\PP_x[T_{(\be,-k)}>t,T_{(\be,1)}>T_{\be}\leq t]-\PP_x[T_{(\be,k)}>t,T_{(\be,1)}>T_{\be}\leq t]\right)\\
S3&=& \sum_{k=1}^{\infty}\left(\PP_x[T_{(\be,-k)}>t,T_{\be}>T_{(\be,1)} \leq t]-\PP_x[T_{(\be,k)}>t,T_{\be}>T_{(\be,1)} \leq t]\right)
\end{eqnarray*}
Then the implication $(T_{\be}\wedge T_{(\be,1)}>t \Rightarrow \forall  k,\; T_{(\be,k)}>t)$ shows that all summands in $S1$ are $0$.
For $S3$ set $a_k=\PP_x[T_{(\be,-k)}>t,T_{\be}>T_{(\be,1)} \leq t]$ and $b_k=\PP_x[T_{(\be,k)}>t,T_{\be}>T_{(\be,1)} \leq t]$. Set $X_u'=X_u\mathbf{1}_{u\le T_{(\be,1)}}+s_{\be,1}X_u\mathbf{1}_{u> T_{(\be,1)}}$. Then $X$ and $X'$ have the same law so
$\label{prem}	a_k=\PP_x[T_{(\be,-k)}'>t,T_{\be}'>T_{(\be,1)}' \leq t].$
For $k\in\Z$, the definition of $X'$ gives
\begin{equation}\label{deux}
T_{(\be,-k)}'=T_{(\be,-k)}\mathbf{1}_{T_{(\be,-k)}\le T_{(\be,1)}} +(T_{(\be,2+k)}\circ\theta_{T_{(\be,1)}}+T_{(\be,1)})\mathbf{1}_{T_{(\be,-k)}> T_{(\be,1)}}\end{equation}
where $\theta$ is the shift operator.
With $k=-1$ this gives $T_{(\be,1)}'=T_{(\be,1)}$. 
With $k=0$ we get $\{T_{\be}'>T_{(\be,1)} \leq t\}=\{T_{\be}>T_{(\be,1)} \leq t\}$ and for all $k$,
\begin{equation}\label{trois}a_k=\PP_x[T_{(\be,-k)}'>t,T_{\be}>T_{(\be,1)} \leq t].\end{equation}
If $k\ge 1$ and $T_{\be}>T_{(\be,1)}$ then $T_{(\be,-k)}\ge T_{\be}>T_{(\be,1)}$, so (\ref{deux}) gives $T_{(\be,-k)}'=T_{(\be,2+k)}\circ\theta_{T_{(\be,1)}}+T_{(\be,1)}$. So (\ref{trois}) becomes
$$a_k=\PP_x[T_{(\be,2+k)}\circ\theta_{T_{(\be,1)}}+T_{(\be,1)}>t,T_{\be}>T_{(\be,1)} \leq t].$$
For $k\ge 0$, $T_{(\be,2+k)}>T_{(\be,1)}$ so $T_{(\be,2+k)}=T_{(\be,2+k)}\circ\theta_{T_{(\be,1)}}+T_{(\be,1)}$ and 
$$a_k=\PP_x[T_{(\be,2+k)}>t,T_{\be}>T_{(\be,1)} \leq t]=b_{2+k}.$$
In this way we get $S3=2\lim_{k\to +\infty} a_k-b_{1}-b_2$. Now $b_{1}=0$, $b_2=a_0$
and since $\{X(s): 0 \leq s \leq t\}$ is almost surely bounded
we have
 $\lim_{k\to +\infty} a_k =\PP_x[T_{\be}>T_{(\be,1)} \leq t]$ so that $$S3=2\PP_x[T_{\be}>T_{(\be,1)} \leq t]-\PP_x[ T_{(\be,1)} \leq t, T_{\be}>t].$$ The same line of reasoning 
gives $S2=0.$ Finally observe that
$$\sum_{k=1}^{\infty}\left(\PP_x[T_{(\be,-k)}>t]-\PP_x[T_{(\be,k)}>t]\right)=
S1+S2+S3.$$ \hfill $\Box$\\

\textbf{Proof of Proposition \ref{ninenine}} Apply Lemmas \ref{sneez} and \ref{cougz} to Theorem \ref{var3z}.

\subsection{Consistency in the different type cases}
\subsubsection{$\widetilde{A}_{k-1}$, $k$ even}
Let us first determine $\I$. If $w_a=\t(d)\s\in W_a^{I}$, then
$$w_a \{ (e_{2i-1}-e_{2i},0),(-e_{2i-1}+e_{2i},-1)\}=$$
$$\{ (e_{\s(2i-1)}-e_{\s(2i)},n),(-e_{\s(2i-1)}+e_{\s(2i)},-1-n)\},$$
where $n=d_{\s(2i-1)}-d_{\s(2i)}$. Thus, $n\le 0$ and $-1-n\le 0$, ie
$n\in \{0,-1\}$. If $n=0$, $d_{\s(2i-1)}=d_{\s(2i)}$ and
$\s(2i-1)<\s(2i)$. If $n=-1$, $d_{\s(2i-1)}=d_{\s(2i)}-1$ and
$\s(2i-1)>\s(2i)$. In any case, $$w_a \{
(e_{2i-1}-e_{2i},0),(-e_{2i-1}+e_{2i},-1)\}=$$
$$\{ (e_{\min(\s(2i-1),\s(2i))}-e_{\max(\s(2i-1),\s(2i))},0),(-e_{\min(\s(2i-1),\s(2i))}+e_{\max(\s(2i-1),\s(2i))},-1)\}.$$
Thus, we identify $\pi=\{\{i_l < j_l\},\, 1\le l\le
p\}\in P_2(k)$ and $A=\{ (e_{i_l}-e_{j_l},0),\,
(-e_{i_l}+e_{j_l},-1)\, ;\, 1\le l\le p\}\in \I$. Then we take $J_a=\{ (e_{2i-1}-e_{2i},0)\, ;\, 1\le i\le
p\}\in \O(\D_a)$. From the previous description of $\I$, (C1) and
(C3) are obvious. Now it is clear that $$U_a=\{\t(d)\s :\ \s \text{ permutes sets } \{1,2\}, \{3,4\},\ldots ,\{k-1,k\}\text{ and }\forall i\le p,\,$$
$$ (d_{\s(2i-1)}=d_{\s(2i)},\, \s(2i-1)<\s(2i)) \text{
or } (d_{\s(2i-1)}=d_{\s(2i)}-1,\,\s(2i-1)>\s(2i))\}.$$
Thus if $\t(d)\s\in U_a$ we can write $\s=\s_1 \s_2$, where $\s_2$
permutes pairs $(1,2),\ldots ,(k-1,k)$ and $\s_1$ is the product of
the transpositions $(\s(2i-1),\s(2i))$ for which
$d_{\s(2i-1)}=d_{\s(2i)}-1$. Then $\varepsilon(\s_2)=1$ from \cite{doc} so
that $\varepsilon(\s)=\varepsilon(\s_1)=(-1)^m$, where $m=|\{i:\
d_{\s(2i-1)}=d_{\s(2i)}-1\}|$. But since $d\in L$, 
\begin{eqnarray}\label{qazal}
0&=&\sum_j d_j=\sum_{i=1}^p \left( d_{\s(2i-1)}+d_{\s(2i)}\right)\\
&=&2 \sum_{i,\, d_{\s(2i-1)}=d_{\s(2i)}} d_{\s(2i)} +2 \sum_{i,\,
d_{\s(2i-1)}=d_{\s(2i)}-1} d_{\s(2i)} -m ,
\end{eqnarray}
which proves that $m$ is even. Hence $\varepsilon(\s_1)=1$. The fact that $\varepsilon_A=(-1)^{c(\pi)}$ comes from the analogous fact in
\cite{doc}.

\textbf{Remark}\label{rema}
In the case of odd $k=2p+1$, the same discussion carries over by adding
singlets to the pair partitions and with $\s(k)=k$ if $\t(d)\s\in
U_a$. But equality (\ref{qazal}) is no longer valid, which explains why
the sign is not well-defined for such $k$.

\subsubsection{The cases $\widetilde{B}_{k}$ and $\widetilde{C}_{k}$}
The argument for the cases $\widetilde{B}_{k}$ and $\widetilde{C}_{k}$ is the same; we give the details in the $\widetilde{B}_{k}$ case.
Let us first suppose $k$ is even, $k=2p$. Suppose $d\in L$, $f$ is a
sign change with support $\bar{f}$ and $\s\in \mathfrak{S}_k$ such that $w_a=\t(d)f\s\in W_a^{I}$. Then,
$$w_a \big\{\, (e_{2i-1}-e_{2i},0),\, (e_{2i},0),\,(-e_{2i-1}-e_{2i},-1)\,\big\}=\big\{ \left( f(e_{\s(2i-1)})-f(e_{\s(2i)}),m-n\right),$$
$$\left(f(e_{\s(2i)}),n\right),\,\left(-f(e_{\s(2i-1)})-f(e_{\s(2i)}),-1-m-n\right)\big\}:=S ,$$
with $m=f(\s(2i-1))d_{\s(2i-1)}$ and $n=f(\s(2i))d_{\s(2i)}$. Thus,
$m-n\le 0,\,n\le 0,\, -1-m-n\le 0$, which forces $m=n=0$ or
$m=-1,n=0$. If $m=n=0$, then $f(e_{\s(2i-1)})-f(e_{\s(2i)})\in\Phi^+$,
$f(e_{\s(2i)})\in\Phi^+$, which implies $\s(2i-1),\s(2i)\notin \f$ and
$\s(2i-1)<\s(2i)$. If $m=-1,n=0$, then $-f(e_{\s(2i-1)})-f(e_{\s(2i)})\in\Phi^+$,
$f(e_{\s(2i)})\in\Phi^+$, which implies $\s(2i-1)\in\f,\, \s(2i)\notin \f$ and
$\s(2i-1)<\s(2i)$. In any case, $$S=\big\{\, (e_{\s(2i-1)}-e_{\s(2i)},0),\,(e_{\s(2i)},0),\,(-e_{\s(2i-1)}-e_{\s(2i)},-1)\big\}$$
and 
\begin{eqnarray*}
W_a^{I}&=\Big\{&\t(d)f\s\in W_a :\ \forall i,\,
\big(\,d_{\s(2i-1)}=d_{\s(2i)}=0,\; \s(2i-1),\s(2i)\notin \f,\\
&&\s(2i-1)<\s(2i)\,\big)
\text{ or } \big(\,d_{\s(2i-1)}=1,\,d_{\s(2i)}=0,\;\s(2i-1)\in\f,\\
&& \s(2i)\notin \f,\;
\s(2i-1)<\s(2i)\,\big)\ \Big\}.
\end{eqnarray*}
Then $\I$ clearly identifies with $P_2(k)$ through the
correspondence between $\pi=\{\{i_l < j_l\},\, 1\le l\le
p\}\in P_2(k)$ and $A=\{ (e_{i_l}-e_{j_l},0),\,(e_{j_l},0),\,
(-e_{i_l}-e_{j_l},-1)\, ;\, 1\le l\le p\}$. So, (C1) and (C3) are
obvious by taking $J_a=\{(e_{2i-1}-e_{2i},0),\,(-e_1-e_2,-1)\}$. Now, 
$$U_a=\{\t(d)f\s\in W_a^{I}: \s \text{ permutes pairs } (1,2),\ldots
,(2p-1,2p)\},$$
so that if $\t(d)f\s\in U_a$,
$\varepsilon(\t(d)f\s)=\varepsilon(f)\varepsilon(\s)=(-1)^{|\f|}$. But $|\f|=\sum_i
d_{\s(2i-1)}=\sum_j d_j$ is even, which proves (C2).\\
For odd $k=2p+1$, $\I$ identifies with $P_2(k)$ through the
correspondence between $\pi=\{\{i_l < j_l\},\, 1\le l\le p;\,\{s\}\}\in P_2(k)$ and $A=\{ (e_{i_l}-e_{j_l},0),\,(e_{j_l},0),\,
(-e_{i_l}-e_{j_l},-1)\, ,\, 1\le l\le p;\,
(e_s,0),\,(-e_s,-1)\}$. Elements $\t(d)f\s\in U_a$ are described in the same
way with the extra condition that $\s(k)=k$ and $d_k=0, k\notin\f$ or
$d_k=1, k\in\f$. So the proof of (C2) carries over.

\subsubsection{The $\widetilde{D}_{k}$ case}
Let us first suppose $k$ is even, $k=2p$. Suppose $d\in L$, $f$ is an even
sign change and $\s\in \mathfrak{S}_k$ such that $w_a=\t(d)f\s\in W_a^{I}$. Then,
\begin{gather*}
w_a \left\{\, (e_{2i-1}-e_{2i},0),\,
(-e_{2i-1}+e_{2i},-1),\,(e_{2i-1}+e_{2i},0)\,(-e_{2i-1}-e_{2i},-1)\,\right\}\\
=\left\{ \left( f(e_{\s(2i-1)})-f(e_{\s(2i)}),m-n\right),\left(-f(e_{\s(2i-1)})+f(e_{\s(2i)}),-1-(m-n)\right),\right.\\
\left.\left(f(e_{\s(2i-1)})+f(e_{\s(2i)}),m+n\right),\left(-f(e_{\s(2i-1)})-f(e_{\s(2i)}),-1-(m+n)\right)\right\}:=S,
\end{gather*}
with $m=f(\s(2i-1))d_{\s(2i-1)}$ and $n=f(\s(2i))d_{\s(2i)}$. Thus
$m-n\le 0,\,-1-(m-n)\le 0,\,m+n\le 0,\,-1-(m+n)\le 0$, which forces $m=n=0$ or
$m=-1,n=0$. If $m=n=0$, then $f(e_{\s(2i-1)})\pm f(e_{\s(2i)})\in\Phi^+$, which implies $\s(2i-1)\notin \f$ and
$\s(2i-1)<\s(2i)$. If $m=-1,n=0$, then $-f(e_{\s(2i-1)})\pm f(e_{\s(2i)})\in\Phi^+$, which implies $\s(2i-1)\in\f$ and
$\s(2i-1)<\s(2i)$. In any case, we have
\begin{gather*}
S=\left\{\,
(e_{\s(2i-1)}-e_{\s(2i)},0),\,(-e_{\s(2i-1)}+e_{\s(2i)},-1)\right.,\\
\left. \qquad (e_{\s(2i-1)})+e_{\s(2i)},0)\,(e_{\s(2i-1)})+e_{\s(2i)},0)\,\right\},
\end{gather*}
and
\begin{eqnarray*}
W_a^{I}&=\Big\{&\t(d)f\s\in W_a :\ \forall i,\,
\big(\,d_{\s(2i-1)}=d_{\s(2i)}=0,\; \s(2i-1)\notin \f,\\
&&\s(2i-1)<\s(2i)\big)\text{ or }
\big(\,d_{\s(2i-1)}=1,\,d_{\s(2i)}=0,\,\s(2i-1)\in\f,\\
&&\,(2i-1)<\s(2i)\,\big)\ \Big\}.
\end{eqnarray*} 
The
correspondence between $\pi=\{\{i_l < j_l\},\, 1\le l\le
p\}\in P_2(k)$ and $A=\{
(e_{i_l}-e_{j_l},0),\,(-e_{i_l}+e_{j_l},-1),\,(e_{i_l}+e_{j_l},0),(-e_{i_l}-e_{j_l},-1)\,
;\, 1\le l\le p\}$ identifies $\I$ with $P_2(k)$. (C1) and (C3) are
obvious with $J_a=\{(e_{2i-1}-e_{2i},0),\, 1\le i\le
p;\,(e_{k-1}+e_k,0)\}$. Moreover, $$U_a=\{\t(d)f\s\in W_a^{I}: \s \text{ permutes pairs } (1,2),\ldots
,(2p-1,2p)\},$$
which makes (C2) easy since $\varepsilon(f)=1$ for $\t(d)f\s\in W_a$.\\
The case of odd $k$ is an obvious modification.

\subsubsection{The $\widetilde{G}_{2}$ case}
Call $\a_1=e_1-e_2$, $\a_2=2e_3-e_1-e_2=\widetilde{\al}$ and take
$J_a=\{(\a_1,0),\,(-\a_2,-1) \}$. We remark that $I$ can be written
\begin{equation}\label{pro}
\{(\a_1,0),\,(-\a_1,-1),\,(\a_2,0),\,(-\a_2,-1)\}\text{ with }
\a_1\text { short}, \a_2 \text{ long,}\ \a_1\perp \a_2.
\end{equation}
If $w_a=\t(d)w\in W_a^{I}$
then $(w\a_i,d)\in\Z$, $(w\a_i,d)\le 0$ and $-1-(w\a_i,d)\le 0$, which imposes
$(w\a_i,d)\in\{0,-1\}$ for $i=1,2$. Thus, $A=w_a I$ can also be written as in
(\ref{pro}) for some $\a_1', \a_2'$. This guarantees condition (C3) and if $J_a\subset A$ then obviously
$\a_1=\a_1'$, $\a_2=\a_2'$ so that $A=I$, which proves condition
(C1). Writing $I$ as in (\ref{pro}) allows us to see that if $w_a=\t(d)w\in
W_a$, then
$w_aI=\{(w\a_1,m_1),(-w\a_1,-1-m_1),(w\a_2,m_2),(-w\a_2,-1-m_2)\}$ where
$m_i=(w\a_i,d)\in\Z$. Since $W$ sends long (short) roots to long (short)
roots, $w_a\in U_a$ implies $w\a_i\in \{\pm \a_i\}$ for $i=1,2$. If
$w\a_i=\a_i$ for $i=1,2$ (respectively $w\a_i=-\a_i$ for $i=1,2$), then
$w=\id$ (respectively $w=-\id$) and $\varepsilon(w)=1$ (recall that $\dim V=2$). If
$w\a_1=\a_1$ and $w\a_2=-\a_2$ then $(\a_1,d)=0$ and $(\a_2,d)=1$. This implies
$d=(-1/6,-1/6,1/3) \notin L$, which is absurd. The same absurdity occurs if
$w\a_1=-\a_1$ and $w\a_2=\a_2$.\\ 

For the determination of $\I$, it is easy to see that the sets of the form
(\ref{pro}) are $I,A_1,A_2$. The sign of the transformation sending
$(\a_1,\a_2)$ to $(e_3-e_1,-2e_2+e_1+e_3)$ is $1$ so that $\varepsilon_{A_1}=-1$ and
$A_2$ is obtained from $A_1$ by transposing $e_1$ and $e_2$, which finishes
the proof.

\subsubsection{The $\widetilde{F}_{4}$ case}
Call $\a_1=e_2-e_3$, $\a_1'=e_3$, $\a_2=e_1-e_4$, $\a_2'=e_4$. Then $I$ can
be written
\begin{equation}\label{prof}
\{(\a_1,0),(-\a_1,-1),(\a_1',0),(\a_2,0),(-\a_2,-1),(\a_2',0)\},
\end{equation}
with $\a_1,\a_2$ long, $\a_1',\a_2'$ short, $\{\a_1,\a_1'\}\perp
\{\a_2,\a_2'\}$ and $(\a_i,\a_i')=-1$. The same kind of reasoning as in the $\widetilde{G}_{2}$ case
shows conditions (C1) and (C3), with $J_a=\{\a_1,\a_2'\}$. Let us prove (C2). If $w_a=\t(d)w\in U_a$, then $w_a I=$
\begin{eqnarray*}
\{(w\a_1,m_1),(-w\a_1,-1-m_1),(w\a_1',m_1'), \qquad\\
(w\a_2,m_2),(-w\a_2,-1-m_2),(w\a_2',m_2')\},
\end{eqnarray*}
with $m_i=(w\a_i,d)$, $m_i'=(w\a_i',d)$. Since $w$ sends long (short) roots to
long (short) roots, necessarily $w\{\a_1',\a_2'\}=\{\a_1',\a_2'\}$ and
$m_1'=m_2'=0$.\\
Suppose
$w\a_i'=\a_i',\, i=1,2$. Since $(w\a_2,\a_1')=(\a_2,\a_1')=0\not= -1$, we have
$w\a_1\in \{\a_1,-\a_1\}$ and $w\a_2\in \{\a_2,-\a_2\}$. If $w\a_1=-\a_1,\,
w\a_2=\a_2$ then $m_1=1$, $m_2=0=m_1'=m_2'$, which leads to $d=(0,1,0,0)\notin
L$, absurd! If $w\a_1=\a_1,\,
w\a_2=-\a_2$, a similar reasoning leads to the absurdity $d=(1,0,0,0)\notin
L$. Hence, $w\a_1=\a_1,\,
w\a_2=\a_2$ or $w\a_1=-\a_1,\, w\a_2=-\a_2$. Then, using the basis
$(\a_1,\a_1',\a_2,\a_2')$, $\varepsilon(w)=1$ is easily checked. 
\\Suppose now $w\a_1'=\a_2',\,w\a_2'=\a_1'$. Similar arguments show that
$w\a_2\in \{\a_1,-\a_1\}$ and $w\a_1\in \{\a_2,-\a_2\}$. If
$w\a_1=\a_2,\,w\a_2=\a_1$ or $w\a_1=-\a_2,\,w\a_2=-\a_1$ then $\varepsilon(w)=1$.
Suppose $w\a_1=\a_2,\,w\a_2=-\a_1$, then $m_1=0,\, m_2=-1$, which, as before,
leads to $d=(0,1,0,0)\notin L$. If $w\a_1=-\a_2,\,w\a_2=\a_1$, then $m_1=-1,\, m_2=0$,
which also gives $d=(1,0,0,0)\notin L$.

\hfill $\Box$\\

\subsection{Proposition \ref{mano}}\label{prpf}
The definition of the Pfaffian is given in the appendix. We refer to \eqref{forma} for even $k$ and \eqref{leen} and for odd $k$.
\hfill $\Box$\\

\subsection{Proposition \ref{librz} 
} \label{expaz}

We will use the following expansions involving the exit time $T_{(0,1)}$ from $(0,1)$ and the hitting times $T_0$ and $T_1$ of 0 and 1 respectively for one-dimensional
Brownian motion: for $(x,t) \in (0,1) \times [0,\infty)$,
\begin{eqnarray}
\phi(x,t)&:=&\P_x(T_{0,1}>t)=\sum_{l\in 2\nn+1} c_l e^{-\la_l t} \sin(\pi l x) \label{fzt}\\
\psi(x,t)&:=&\P_x(T_{0,1}>t)+2\P_x(T_0>T_1)-2\P_x(T_0>T_1>t)\nonumber\\&=&
\sum_{l\in 2\nn} c_l e^{-\la_l t} \sin(\pi l x) \nonumber
\end{eqnarray}
with 
$c_l=4/(l\pi)$, $\la_l=(l\pi)^2/2$ and the definition 
$c_l \sin (\pi l x) = 2 x$ when $l=0$. The first expansion may be found in, for example, \cite{MR1329542}; in the case of even $k$, it may be used to rewrite (\ref{forma}) into the form \eqref{markus}. The second expansion is obtained using 

\nw{sor}[y1]{Lemma}
\bg{sor} If $X$ is Brownian motion and $\be=e_i-e_j$ then 
\label{sorz}$$\PP_x[T_{\be}> T_{(\be,1)}>t]=2\sum_{n=1}^{\infty} \frac{(-1)^{n+1}}{\pi n}
e^{- \pi^2 n^2 t}\sin (\pi n x_{ij})$$
where $x_{ij}=x_i-x_j \in (0,1)$.
\end{sor}
\textbf{Proof} 
The series satisfies the diffusion equation for $(x_{ij},t)\in (0,1)\times (0,\infty)$, takes the value 0 if $x_{ij} \in \{0,1\}$, and equals $x_{ij}$ if $t=0$.
$X_{ij}:=X_i-X_j$ is a Brownian motion with the same diffusion coefficient.  
Therefore by applying for example Theorem 4.14 of \cite{MR1329542},
the series equals 
\begin{eqnarray*}
\EE_x[X_{ij}(t); T_{\be}\wedge T_{(\be,1)}>t]&=&\EE_x[\PP_x[T_{\be}> T_{(\be,1)}>t
|X(t),\mathbf{1}_{T_{\be}\wedge T_{(\be,1)}>t}]]\\
&=&\PP_x[T_{\be}> T_{(\be,1)}>t].
\end{eqnarray*} \hfill $\Box$\\
We record the following corollary, which follows from integration, interchanging integration with summation, and inversion of Fourier series:
\nw{cc}[y1]{Corollary}
\bg{cc}Under the conditions of Lemma \ref{sorz}, 
\label{ccz}
$$\int_0^{\infty}\PP_x[T_{\be}> T_{(\be,1)}>t]
dt=\frac{1}{6}x_{ij}(1-x_{ij}^2).$$
\end{cc}

In the case of odd $k$, the second expansion in \eqref{fzt} may be 
inserted in Proposition \ref{ninenine} to give \eqref{markus}: 
\begin{eqnarray}
\P_x(\tT>t)=\sum_{\pi=\{\{i_s<j_s\},\,1\le s\le m\} } (-1)^{c(\pi)}
\prod_{s=1}^{m} \left(\sum_{l\in \mathbb{O}} c_l e^{-2\la_l t} \sin(\pi l x_{i_s j_s}) \right)\label{markus}\\ 
\quad= \sum_{\pi=\{\{i_s<j_s\},\,1\le s\le m\}} (-1)^{c(\pi)}
\sum_{l\in \mathbb{O}^m} e^{- \pi^2 (l_1^2+\cdots+l_m^2) t} \prod_{s=1}^{m} c_{l_s} \sin(\pi l_s x_{i_s j_s})\nonumber
\end{eqnarray}
 for $x \in \cA$, where $m=\lfloor k/2 \rfloor \in \N$, $x_{ij}=x_i-x_j$, $\mathbb{O}=2\N+1$ if $k$ is even and
$\mathbb{O}=2\N$ if $k$ is odd.
Now for $\pi=\{\{i_s<j_s\},\,1\le s\le m\}$ define 
\begin{eqnarray}
G_r(x,\pi)=\sum_{l\in \mathbb{O}^m,\,N(l)=r}\; \prod_{s=1}^{m} c_{l_s} \sin(\pi l_s
x_{i_s j_s}) \label{germ}
\end{eqnarray}
where $N(l)=l_1^2+\cdots+l_m^2$, and let
$F_r(x)=\sum_{\pi\in P_2(k)} (-1)^{c(\pi)} G_r(x,\pi).$
(Since the sum defining $G_r(x,\pi)$ runs over a $\mathfrak{S}_m$-invariant set of indices, it does not depend on
the enumeration of the blocks of $\pi$ but only on $\pi$ itself.) With those definitions we can write
\begin{equation}\label{alc_expansion}
\P_x(\tT>t)=\sum_{r>0} e^{-\pi^2 rt} F_r(x)
\end{equation}
(note that by Proposition 2.4 of \cite{doc}, 
$\sum_{\pi \in P_2(k)}(-1)^{c(\pi)}\prod_{s=1}^m x_{i_s j_s}=0$
and so the terms corresponding to $r=0$ cancel.)
As for expectations, we have
\begin{eqnarray}
\mathbb{E}_x(\tT)=\int_0^\infty \P_x(\tT>t)\,dt
=\sum_{r>0} \frac{1}{r\pi^2} F_r(x)\nonumber
\end{eqnarray}
and the result follows. \hfill $\Box$\\

When $k=2$ the previous formula becomes
\begin{eqnarray}
\mathbb{E}_x(\tT)=\sum_{n\in\N}\frac{4}{\pi^{3}} \frac{\sin\left(\pi(2n+1)x_{12}\right)}{(2n+1)^3}=\frac{1}{2} x_{12}(1-x_{12}), \label{timez}
\end{eqnarray}
$0<x_{12}<1,$
which is a well-known formula in Fourier series. 
When $k=3$ we may use the above and Corollary \ref{ccz} to obtain
\begin{eqnarray}
\mathbb{E}_x(\tT)=\sum_{\pi=\{i_s<j_s\}}(-1)^{c(\pi)}\sum_{n\in\N}\frac{4}{\pi^{3}} \frac{\sin\left(2\pi nx_{ij}\right)}{(2n)^3}=x_{12}x_{23}(1-x_{13}),\label{alwayz}
\end{eqnarray}
$0<x_{ij}<1.$ It is easy to check that \eqref{timez} and \eqref{alwayz} both solve Poisson's equation
$\frac{1}{2}\Delta u=-1$
inside the interval and an equilateral triangle respectively and vanish on the boundary, which confirms that they are the expected exit times for Brownian motion from these domains. Formula \eqref{alwayz} has also been obtained using scaling limits for random walks (see \cite{MR2023644, MR1975514} ).

\subsection{The reflection principle and De Bruijn Formulae}

\subsubsection{Proposition \ref{certz}} 
From \eqref{barf},
if $T_A$ is the exit time of Brownian motion from $F_A$ then
\begin{eqnarray}
\PP_x[T_A>t]=\int_{F_A} \sum_{\om \in W_A}\ep(\om)p_t(x,\om y)dy \label{finl}
\end{eqnarray}
where $p_t$ is the Brownian transition density and $x \in F_A$. The finite case was proved in \cite{doc}; in the affine case it is easy to check that the same proof applies.

\subsubsection{Propositions \ref{greenz} and \ref{green2z}} 
We treat first the case of odd $k$.
Let $\be \in \Phi^+=\{e_i-e_j:1 \leq i<j \leq k\}$ and $x \in \cA$. Then $(x,\be) \in (0,1)$ and from~\eqref{finl}, for $k \geq 1$
\begin{eqnarray*}
\PP_x[T_{(\be,k)}>t]&=&\int_{\<y,\be\ra<k}p_t(x,y)-p_t(x,s_{\be}y+k\be)dy\\
&=&\int_{\<u,\be\ra>-k}p_t(x,s_{\be}u)-p_t(x,u+k\be)du
\end{eqnarray*}
where $u=s_{\be}y$. Also if $k \leq 0$ then
\begin{equation}
\PP_x[T_{(\be,k)}>t]=\int_{\<y,\be\ra>k}p_t(x,y)-p_t(x,s_{\be}y+k\be)dy. \label{tsho}
\end{equation}
Write $\be=e_i-e_j$. Rewriting Theorem \ref{var3z} using Lemma \ref{sneez}, equation \eqref{tsho} and using the identification of $\ii$ with $P_2(k)$ from section \ref{drop} we have 
\begin{eqnarray}
\PP_x[\tT>t]&=&\sum_{\pi \in P_2(k)}(-1)^{c(\pi)} \prod_{\{i<j\} \in \pi}
\Bigg( \int_{y_i>y_j} p_{ij}(0) dy_i dy_j 
\nonumber\\&& \quad \label{ooz}
+\sum_{k=1}^{\infty}\int_{y_i-y_j>-k} p_{ij}(0)+p_{ij}(k) dy_i dy_j\Bigg)
\end{eqnarray}
where $p_{ij}(k)=\psi(x_i,y_i+k)\psi(x_j,y_j-k)-
\psi(x_i,y_j-k)\psi(x_j,y_i+k)$ and $\psi(x,y)=\frac{1}{\sqrt{2\pi t}}e^{-(x-y)^2/2t}$.
Now $\int_{-k<y_i-y_j<k}p_{ij}(0)dy_i dy_j=0$ and 
making the substitution $(u_i,u_j)=(y_i+k,y_j-k)$ we have
$$\int_{y_i-y_j>-k}p_{ij}(k)dy_i dy_j=
\int_{u_i-u_j>k}p_{ij}(0)du_i du_j,$$ so the infinite sum in \eqref{ooz} may be written
$2 \sum_{k=1}^{\infty}\int_{y_i-y_j>k}p_{ij}(0) dy_i dy_j.
$
From~\eqref{finl} we have the alternative expression
$$\PP_x[\tT>t]=\int_{\cA} \sum_{\om \in W_a}\ep(\om)p_t(x,\om y)dy
$$
so integrating both expressions over $\rr^k$ with respect to $f_i(x_i)dx_i$, $i=1,\ldots , k$ and applying Fubini's theorem,
\begin{eqnarray*}
&&\int_{\cA} \sum_{\om \in W_a}\ep(\om)P_tf(\om y)dy=\sum_{\pi \in P_2(k)}(-1)^{c(\pi)} 
\\
&& \;
\int_{\rr}f_{l_{\pi}}\prod_{\{i<j\} \in \pi}
\Bigg( \int_{y_i>y_j} P_{ij}
dy_i dy_j 
+2 \sum_{k=1}^{\infty}\int_{y_i-y_j>k} P_{ij}
dy_i dy_j
\Bigg)
\end{eqnarray*}
where $\{l_{\pi}\}$ is the singlet in the partition $\pi$ and
$P_{ij}=P_tf_i(y_i)P_tf_j(y_j)-
P_tf_i(y_j)P_tf_j(y_i)$.
To complete the proof for the case of odd $k$ we obtain 
uniform bounds in $t$ to justify the use of dominated convergence to let $t \to 0$ inside the infinite sum, and finally apply the definition of the Pfaffian.
Dividing the domain of integration into $(-k+\sqrt t, k-\sqrt t)$ and its complement and applying the bound $\int p_{ij}(0) dy_i dy_j \leq 2$ on the latter we have for $t < 1/4$
\begin{eqnarray*}
&&\left|\int_{y_i-y_j> k}P_{ij} dy_idy_j\right|\leq
\int_{x_i,x_j \in \rr} \int_{y_i-y_j>k}|p_{ij}(0)f_i(x_i)f_j(x_j)|
dy_i dy_j dx_i dx_j \\&& \;
\leq \int_{x_i-x_j<k-\sqrt t} 
\int_{y_i-y_j>k} \psi(x_i,y_i)\psi(x_j,y_j)(|f_i(x_i)f_j(x_j)|+|f_i(x_j)f_j(x_i)|)
\\ &&\quad  dy_i dy_j dx_i dx_j+2
\int_{x_i-x_j \in (-\infty,-k+1/2)\cup(k-1/2,\infty)} |f_i(x_i)f_j(x_j)| dx_i dx_j
\end{eqnarray*}
and $\int_{x_i-x_j \in (-\infty,-k)\cup(k,\infty)}
|f_i(x_i)f_j(x_j)| dx$ is summable in $k$ by assumption.
The standard estimate for the tail of the Gaussian distribution gives $$\int_{y_i-y_j>k}\psi(x_i,y_i)\psi(x_j,y_j)dy_i dy_j\leq
e^{-(k-(x_i-x_j))^2}$$ when $x_i-x_j<k-\sqrt t$, and 
$\int_{x \in \mathbb{R}} e^{-(k-(x_i-x_j))^2}|f(x)| dx$ is summable in $k$.

When $k$ is even we have a consistent subset $I$ as described in section \ref{drop} and so Proposition \ref{certz} applies. The proof is similar to that in section 7.6.1 of \cite{doc}, with the difference that here we have the bijection
$$(l \in L_{\pi},\eta \in\{\pm 1\}^{\pi}) \mapsto w_{l,\eta}=\tau(l) \prod_{\{i<j\}\in \pi}
\tau_{ij}^{\eta_{ij}'} \in W_A$$
where $\pi \in P_2(k)$ is the pair partition associated with $A \in \ii$, and $L_\pi$ is the coroot lattice associated with the affine Weyl group $W_A$; and now
$F_A$ corresponds with 
$F_{\pi}=\cap_{\{i<j\} \in \pi}\{y:0<y_i-y_j<1\}.$
\hfill $\Box$\\

\section{Eigenfunctions and eigenvalues for alcoves}
\label{eigz}
It follows from equation (\ref{alc_expansion}) that $F_r$ is a real eigenfunction for the Dirichlet Laplacian on the alcove of type $\tA_{k-1}$, with eigenvalue $-2\pi^2 r$.
As an example, when $k=3$ the alcove is the equilateral triangle and we have 
\begin{eqnarray*}
\frac{1}{c_{2n}} F_r(x)=\left\{ \begin{array}{ll}\sin (2\pi n x_{12})+\sin (2\pi n x_{23})-\sin (2\pi n x_{13}) & \mbox{ if } r=4n^2\\0 & \mbox{ otherwise, }\end{array} \right.
\end{eqnarray*}
giving the eigenfunctions with simple eigenvalues (see \cite{MR586910}), a feature which can be anticipated from the symmetry of the equilateral triangle. B\'erard \cite{MR570879} obtained a general formula for the eigenfunctions of  the Dirichlet and Neumann Laplacians for alcoves of any type, and we provide here a characterisation of the real eigenfunctions.

Defining
\begin{eqnarray} \label{efns}
	f_p(x) = \sum_{w \in W}\ep(w)\exp(2\pi i  \left< x, wp\right>), \;
	g_p(x)  = \sum_{w \in W}\exp(2\pi i  \left< x, wp\right>),
\end{eqnarray}
the eigenfunctions for the Dirichlet Laplacian on $\cA$ are $\{f_p: p \in \mathcal{P} \cap \cC\}$, where $\ve(w)=\det w$ and $\mathcal{P}=\{x \in V: \left< \a^{\vee},x \right> \in \Z \; \forall \; \a \in \Phi\}$, and the eigenfunctions for the Neumann Laplacian on $\cA$ are $\{g_p: p \in \mathcal{P} \cap \overline{\cC}\}$. 

\textbf{Remark} It is immediate from \eqref{efns} that if $g_{p}$ is real then for every $y \in \cA$ we have $g_{p}(y)<\sup_{x \in \partial \cA}g_{p}(x)$. The `Hot Spots' conjecture of J. Rauch (see \cite{MR1694534}) is therefore true for alcoves. Note that in the two-dimensional case, the alcoves are the equilateral triangle and the right triangles with an angle of either $\pi/4$ or $\pi/3$.

\newtheorem{comp}[y1]{Proposition}
\begin{comp}
(i) For $p \in \mathcal{P} \cap \cC$, the eigenfunction $f_p$ of the Dirichlet Laplacian on $\cA$ is real iff 
\begin{equation}
\label{exiw}
\exists \; w_1 \in W \mbox{ such that } w_1p=-p.
\end{equation} If \eqref{exiw} holds then, up to a constant factor,\begin{equation*}
f_p(x)=\sum_{\substack{w \in W}}
\ve(w) \mbox{cs} (2\pi \left<x,wp\right>)
\end{equation*}
where $cs=\sin$ if $\ve(w_1)=-1$ and $cs=\cos$ if $\ve(w_1)=1$. 

(ii) For $p \in \mathcal{P} \cap \overline{\cC}$, the eigenfunction $g_p$ of the Neumann Laplacian on $\cA$ is real iff \eqref{exiw} holds and then, up to a constant factor, \begin{equation*}
g_p(x)=\sum_{\substack{w \in W}}
 \cos 2\pi \left<x,wp\right>.
\end{equation*}
\end{comp} 
\proof
(i) We have \begin{eqnarray*}
f_p(x) & = &  \sum_{w \in W}\ve(w)\cos 2\pi \left< x, wp\right>+i\sum_{w \in W}\ve(w) \sin2\pi \left< x, wp\right>.
\end{eqnarray*}
Suppose first that $w_1p=-p$ for some $w_1 \in W.$ Then by conjugation, for any $w \in W$ there exists $v_w \in W$ such that $v_w(wp)=-wp$. The orbit $Wp$ may therefore be partitioned into pairs $\{wp,-wp\}$, and 
\begin{equation*}
\mbox{cs} (2\pi \left<x,wp\right>) \pm \mbox{cs}(2\pi \left<x,-wp\right>)=0
\end{equation*}
where $\pm=+,-$ if cs $=\sin, \cos$ respectively. The sufficiency of condition \eqref{exiw} is proved by noting that $ \forall \; w \in W, \;\ve(v_ww)=\ve(v_w)\ve(w)=\ve(w_1)\ve(w)$.

Conversely, suppose that 
\begin{equation}
\label{impart}
 \sum_{w \in W}\ve(w)cs 2\pi \left< x, wp\right> = 0 \quad \forall \; x \in V.
\end{equation}
By restricting $x$ to a ray $x=tr \, (t \in \R)$ chosen such that $\left< r, wp\right>=\left< r, p\right>$ only when $w=\text{Id}$ and $\left< r, wp\right>=-\left< r, p\right>$ only when $wp=-p$ and appealing to linear independence, we conclude that \eqref{exiw} holds.
Part (ii) is proved similarly.\hfill $\Box$

Using standard facts about the longest element of a Weyl group (see \cite{MR1838580}) we obtain
\newtheorem{class}[y1]{Corollary}
\begin{class}\label{cbc}
For the cases $W=A_1,B_k, C_{k}, D_{2k},E_7,E_8,F_4,G_2,H_3$ and $H_4$, all the eigenfunctions of the Laplacian on $\cA$ with Dirichlet or Neumann boundary conditions are real. In all other cases, the eigenfunctions $f_p,g_p$ given by \eqref{efns} are real iff $p=\tau(p)$, where $\tau$ is the unique involution of the Coxeter graph of $W$.
\end{class}
The root systems covered by the second case of Corollary \ref{cbc} are

\textbf{Type $A_{k-1}$, $k > 2$}
Here $\tau(e_i-e_{i+1})=e_{k-i}-e_{k-i+1}$ and so we require $p=\sum_{i=1}^{k-1} a_i(e_i-e_{i+1})$ with $a_i=a_{k-i} \; \forall \; 1 \leq i \leq k-1$. 

\textbf{Type $D_{2k+1}$} Here $\tau$ leaves $e_i-e_{i+1}$ invariant for $1 \leq i \leq 2k-1$, and $\tau(e_{2k}-e_{2k+1})=e_{2k}+e_{2k+1}$. We therefore require $p=\sum_{i=1}^{2k} a_i(e_i-e_{i+1})+a_{2k+1}(e_{2k}+e_{2k+1})$ with $a_{2k}=a_{2k+1}.$

\textbf{Proof of Proposition \ref{trigoharmoyou}}
Defining $\rho=\frac 12 \sum_{\al\in\Phi^+} \al$, we have $\rho \in \mathcal{P} \cap \cC$ (see for example \cite{MR1838580}). Then setting $p=\rho$ in \eqref{efns},
the Weyl identity gives that 
up to a constant factor,
\begin{align}
	f_{\rho}(x)=\prod_{\a \in \Phi^+}\sin 
	( \pi \langle 
	x,\a \rangle).
\end{align}
The next lemma establishes the final claim of Proposition \ref{trigoharmoyou}.
\nw{yan7}[y1]{Lemma}
\begin{yan7}\label{trigolem}
Suppose that $F(X)=F(X_j)_{j \in J}$ is a polynomial in the $(\sin X_j, \cos X_j)_{j \in J}$ which vanishes whenever $\sin X_j$ vanishes.
Then $\sin X_j$ divides $F(X)$ in the ring of trigonometric polynomials.
\end{yan7}
\proof Let $F(X)=P\left(e^{iX_j},e^{-iX_j}
\right)_{j \in J} \in R:=\C[e^{iX_j},e^{-iX_j}\, ;\,j \in J]$. The given cancellation property assures that $P$ is divisible in $R$ by the monic polynomial $e^{iX_{j}}-1$ and that the quotient is divisible by
$e^{iX_{j}}+1$. Hence $P$ is divisible in $R$ by $\frac {e^{2iX_{j}}-1}{2ie^{iX_{j}}}=\sin X_j$.  \hfill $\Box$ 

Since the eigenfunctions are alternating under the action of the affine Weyl group (see for example \cite{MR570879}), and putting $J=\Phi^+$ and $X_{\a}=\pi \langle \a,x\rangle$, the Lemma applies. Using continuity and Lemma \ref{trigolem} again establishes the final claim of Proposition \ref{trigoharmoyou}. 
\hfill $\Box$

\textbf{Remark} In the type $\tA$ case, the principal eigenfunction was obtained by Hobson and Werner in \cite{MR1405497}; we give a direct proof in the appendix. See also \cite{demni07}.

\section{Appendix}

\subsection{Direct proof of Proposition \ref{trigoharmoyou} in the type $\tA$ case}
Set $x_{ij}=x_i-x_j$ and $h(x)=\prod_{1\le i<j\le k} \sin x_{ij}$. Computation of the logarithmic derivative gives
$$\partial_i h=h \sum_{j\,(\not= i)} \frac{\cos x_{ij}}{\sin x_{ij}},$$
which yields
\begin{eqnarray*}
\partial_i^2 h&=&h \left\{ \sum_{j,\,l\,(\not= i)} \frac{\cos x_{ij} \cos x_{il}}{\sin x_{ij} \sin x_{il}} +\sum_{j\,(\not= i)} \left(-1-\frac{\cos^2 x_{ij}}{\sin^2 x_{ij}}\right)\right\}\\
&=&h \left\{ \sum_{j\not= l\,(\not= i)} \frac{\cos x_{ij} \cos x_{il}}{\sin x_{ij} \sin x_{il}} -(k-1)\right\},\end{eqnarray*}
so that $\De h=h (S(x)-k(k-1))$ with 
$$S(x)=\sum\nolimits' \frac{\cos x_{ij} \cos x_{il}}{\sin x_{ij} \sin x_{il}},$$
where $\sum\nolimits'$ runs over $i,j,l$ pairwise distinct. By circular permutation, we get 
\begin{eqnarray*}
&&3S(x)= \sum\nolimits' \frac{\cos x_{ij} \cos x_{il}}{\sin x_{ij} \sin x_{il}}+\frac{\cos x_{jl} \cos x_{ji}}{\sin x_{jl} \sin x_{ji}}+\frac{\cos x_{li} \cos x_{lj}}{\sin x_{li} \sin x_{lj}}\\
&&\; =\sum\nolimits' \frac{\cos x_{ij} \cos x_{il} \sin x_{jl} -\cos x_{jl} \cos x_{ij} \sin x_{il}+ \sin x_{ij}\cos x_{il} \cos x_{jl}}{\sin x_{ij} \sin x_{il}\sin x_{jl}}.
\end{eqnarray*}
But trigonometry shows that each term in the previous sum equals $-1$, so that $S(x)=-k(k-1)(k-2)/3$, which concludes the proof.

\subsection{The Pfaffian}

For completeness we define the Pfaffian. 
If car $\K \neq 2$, any skew-symmetric matrix $A \in \cM_{n} (\K)$ can be written $A = PDP^{t} $ with $P \in GL(n, \K)$, 
$D = \text{diag}(B_{1} , \ldots , B_{q} )$ and 
$B_{l} = 0 \in \K$ or $B_{l}= J = (j- i )_{1 \leq i,j \leq 2} \in \cM_{2} (\K)$. Hence, if n is odd, $\det A = 0$. If $n$ is even, one can use the previous decomposition to prove 

\nw{pfuni}[y1]{Proposition} 
\bg{pfuni} There exists a unique polynomial $\pf \in \Z [X_{ij}, 1 
\leq i < j \leq n]$ 
such that if $A = (a_{ij})$ is a skew-symmetric matrix of size $n$, $\det A = \pf (A)^{2}$ 
and $\pf ({\rm diag}(J , . . . , J )) = 1$.\end{pfuni}
 The Pfaffian has an explicit expansion in terms of the matrix coefficients:
 
\nw{pfdef}[y1]{Proposition}
\bg{pfdef}$$\pf (A) = \sum_{\pi \in P_{2}(n)}(-1)^{c(\pi )} \prod_{\{i <j\} \in \pi}
a_{ij} = \frac{1}{2^{n} (n/2)!} \sum_{\s \in \gS_{n}}
\ep(\sigma)
\prod_{i=1}^{n-1}a_{\sigma (i) \sigma (i+1)}.$$
\end{pfdef}
For more on Pfaffians and their properties, see \cite{MR1606831, MR1069389}. 

\textbf{Acknowledgements.} The authors would like to thank Neil O'Connell for several stimulating conversations and Rodrigo Ba\~nuelos for pointing out the connection to the `Hot Spots' conjecture. They also thank an anonymous referee for careful reading of the manuscript, which led to improved presentation. Research supported in part by Science Foundation Ireland, Grant Number SFI 04/RP1/I512.

\bibliography{alcove}
\end{document}